\documentclass[reqno,12pt]{amsart} 
\usepackage{hyperref}
\usepackage{comment}

\newtheorem{theorem}{Theorem}[section]
\newtheorem{proposition}[theorem]{Proposition}
\newtheorem{corollary}[theorem]{Corollary}

\theoremstyle{remark}
\newtheorem{remark}{Remark}[section]

\newtheorem{example}{Example}[section]

\numberwithin{equation}{section}

\allowdisplaybreaks[4]

\author{Michael J.\ Schlosser}
\address{Fakult\"at f\"ur Mathematik, Universit\"at Wien,
Oskar-Morgenstern-Platz~1, A-1090 Vienna, Austria}
\email{michael.schlosser@univie.ac.at}
\urladdr{http://www.mat.univie.ac.at/{\textasciitilde}schlosse}
\thanks{The author's research was partly supported by
FWF Austrian Science Fund grant 
  \href{https://doi.org/10.55776/P32305}{10.55776/P32305}.}

\title{Bilateral $q$-ultraspherical functions}



\newcommand{\ta}{\theta}

\newcommand{\ba}{\beta}
\newcommand{\ga}{\gamma}
\newcommand{\da}{\delta}

\newcommand{\rd}{\,\mathrm d}
\newcommand{\ri}{\mathrm i}
\newcommand{\C}{\mathbb C}

\newcommand{\Z}{\mathbb Z}

\begin{document}

\begin{abstract}
We introduce the bilateral $q$-ultraspherical functions, a bilateral-series
extension of the continuous $q$-ultraspherical polynomials.  They are defined
by specific bilateral basic hypergeometric ${}_2\psi_2$ series, are analytic in
the variable $x=\cos\theta$, and depend on two parameters $\beta$ and $\gamma$
and on a base $q$.  We derive a product formula for their bilateral generating
function, a three-term recurrence relation, their transformation under the
Askey--Wilson divided difference operator, three weight-based Rodrigues-type
formulae, and explicit large-order asymptotic expansions.  The main results are
full orthogonality relations with respect to explicit orthogonality functionals
involving analytic mass aggregates.  We also obtain shifted orthogonality
relations and a bilateral Chen--Liu type mixed orthogonality formula.  In the
limit $\gamma\to1$, the construction and identities reduce to the classical
results for the continuous $q$-ultraspherical polynomials.
\end{abstract}

\dedicatory{Dedicated to the memory of Richard A.\ Askey}

\maketitle

\section{Introduction}\label{secintro}
Orthogonal and $q$-orthogonal polynomials are central objects in the theory of
special functions; for the one-variable theory, see, for instance, \cite{Sz} and
\cite{Ibook}.  Among the most important examples are the \textit{classical}
orthogonal and $q$-orthogonal families, namely those of hypergeometric and
basic hypergeometric type appearing in the Askey and $q$-Askey schemes
\cite{KLS,KS}. A striking feature of these families is their close connection
with explicit identities for (basic) hypergeometric series. Many summation and
transformation formulae are responsible for various fundamental properties that
orthogonal polynomials possess. It is therefore natural to reverse the point
of view and ask whether a given (basic) hypergeometric identity gives rise to
useful orthogonal polynomials, or more generally to useful orthogonal functions.

Richard Askey had a deep interest in orthogonal and $q$-orthogonal polynomials
and made many fundamental contributions to the subject.  His work with James
Wilson led to the Askey--Wilson polynomials \cite{AW}, the top family in the
$q$-Askey scheme.  Another important contribution of his, joint with Mourad Ismail,
concerns the continuous $q$-ultraspherical polynomials \cite{AI}, the
$q$-analogues of the Gegenbauer polynomials.  These polynomials had been
introduced by Rogers~\cite{R3} in the late nineteenth century; he derived
several of their remarkable properties but apparently was not aware of their
orthogonality.  Askey and Ismail proved their orthogonality with respect to a
positive measure and developed many further properties of this very
important family of special functions.

The continuous $q$-ultraspherical polynomials are highly relevant to
symmetric functions.  In particular, the Macdonald polynomials \cite[Ch.~VI]{Mac},
which are symmetric functions in variables $x_1,\ldots,x_r$ with coefficients
rational in $q$ and $t$, reduce in the two-variable case to the continuous
$q$-ultraspherical polynomials.  Likewise, the Pieri formula for Macdonald
polynomials \cite[p.~331]{Mac} reduces in this case to Rogers' linearization
formula (see \eqref{linear}). It is clear that a solid knowledge about
the continuous $q$-ultraspherical polynomials assists in the study 
of various extensions, including multivariable ones.

The purpose of this paper is to study a \textit{bilateral series extension} of
the continuous $q$-ultraspherical polynomials.  We call the resulting objects
\textit{bilateral $q$-ultraspherical functions} and denote them by
$C_n(x;\ba,\ga\,|\,q)$.  They are defined as explicit multiples of
${}_2\psi_2$ series, contain an additional parameter $\ga$, and specialize to
the ordinary continuous $q$-ultraspherical polynomials when $\ga\to1$.  We show
that many of the classical structural properties survive in this bilateral
setting.  In particular, the bilateral generating function has a closed product
form, from which we derive a three-term recurrence relation and a simple
formula for the action of the Askey--Wilson divided difference operator, as
well as three bilateral Rodrigues-type formulae involving the first
two-factor weight, its dual two-factor weight, the four-factor weight, and
the corresponding shifted seeds $C_0$.
We also prove mixed and quasi-linearization formulae and obtain asymptotic
expansions of $C_n(x;\ba,\ga\,|\,q)$ as $n\to\infty$ and
$n\to-\infty$.  The central result of the paper is the full orthogonality
of the bilateral $q$-ultraspherical functions with respect to explicit
orthogonality functionals involving analytic mass aggregates.  As an
application of the full orthogonality functional, we also derive a bilateral
analogue of Chen--Liu's mixed orthogonality formula for two different
continuous $q$-ultraspherical families~\cite{CL}.  All these results reduce to
the corresponding classical statements for the continuous $q$-ultraspherical
polynomials in the limit $\ga\to1$.

The paper is organized as follows.  Section~\ref{sec:pre} recalls the basic
hypergeometric material needed later, especially identities for $_2\psi_2$
series.  Section~\ref{sec:c} reviews the continuous $q$-ultraspherical
polynomials.  In Section~\ref{sec:bilf} we introduce their bilateral extension
and derive its basic properties, including the generating functions, recurrence
relation, Askey--Wilson operator action, three Rodrigues formulae, and a special
integral.  Section~\ref{sec:asym} gives asymptotic expansions valid for
arbitrary order, and Section~\ref{sec:shifted} records two
shifted orthogonality relations.  Section~\ref{sec:orth} contains the main
orthogonality theorems: full orthogonality for the bilateral
$q$-ultraspherical functions with respect to explicit orthogonality functionals
involving analytic mass aggregates.  Their proof uses the terminating
specializations $\ga=q^s$, the ordinary continuous $q$-ultraspherical
orthogonality relation together with its mass points, and analytic continuation
in $\ga$.  Section~\ref{sec:linear} proves mixed and
quasi-linearization formulae and records their consequences for the residue
mass aggregates.  Section~\ref{sec:chenliu} establishes a
bilateral Chen--Liu type mixed orthogonality formula by reducing first to the
terminating lattice and then continuing meromorphically.  Finally,
Section~\ref{sec:outlook} outlines a possible multilateral Macdonald-type
extension, while Appendix~\ref{app:numerics} discusses numerical checks of the
orthogonality relations.
 
Concerning our notation for fractions, we follow the convention used
in the textbook by Gasper and Rahman \cite{GR} that
all factors that appear after the slash symbol `$/$' are understood
to be part of the denominator.  For instance, `$dq/abz$' (appearing
in \eqref{22tgl}) is meant to stand for `$dq/(abz)$', etc.

\section{Preliminaries on basic hypergeometric series}\label{sec:pre}

\subsection{Notation and convergence}

Let $\Z$ denote the set of integers.  Throughout, we fix $q$ with
$0<|q|<1$.  Occasionally, when positive definiteness is relevant, we further
restrict $q$ to be real and positive, $0<q<1$.  We refer to $q$ as the
``base''.
For a parameter $a\in\C$ and $k\in\Z$, the
$q$-shifted factorial is defined by
\begin{equation*}
(a;q)_k:=\frac{(a;q)_\infty}{(aq^k;q)_\infty},\\
\qquad\text{where}\qquad
(a;q)_\infty:=\prod_{j\ge 0}(1-aq^j).
\end{equation*}
For brevity, we frequently use the following shorthand conventions:
\begin{align*}
(te^{\pm \ri\ta};q)_\infty&=(te^{\ri\ta},te^{-\ri\ta};q)_\infty,\\*
(a_1,\dots,a_m;q)_k&=(a_1;q)_k\dots(a_m;q)_k,\qquad\qquad
k\in\Z\cup\{\infty\}.
\end{align*}

Following Gasper and Rahman~\cite{GR}, basic hypergeometric
${}_r\phi_s$ series with $r$ upper parameters
$a_1,\dots,a_r$, $s$ lower parameters $b_1,\dots,b_s$, base $q$, and
argument $z$ are defined by
\begin{equation}
{}_r\phi_s\!\left[\begin{matrix}
a_1,\dots,a_r\\b_1,\dots,b_s
\end{matrix};q,z\right]:=\sum_{k=0}^\infty
\frac{(a_1,\dots,a_r;q)_k}{(q,b_1,\dots,b_s;q)_k}
\left((-1)^kq^{\binom k2}\right)^{1+s-r}z^k.
\end{equation}
Such a series terminates if one of its upper parameters, say,
$a_r$, is of the form $q^{-n}$, where $n$ is a nonnegative integer.
If the series does not terminate, it converges for $r<s+1$ and, in the
borderline case $r=s+1$, for $|z|<1$.

Bilateral basic hypergeometric ${}_r\psi_s$ series are defined by
\begin{equation}
{}_r\psi_s\!\left[\begin{matrix}
a_1,\dots,a_r\\b_1,\dots,b_s
\end{matrix};q,z\right]:=\sum_{k=-\infty}^\infty
\frac{(a_1,\dots,a_r;q)_k}{(b_1,\dots,b_s;q)_k}
\left((-1)^kq^{\binom k2}\right)^{s-r}z^k.
\end{equation}
Such a series terminates from above if one of its upper parameters,
say, $a_r$, is of the form $q^{-n}$,
and it terminates from below if one of its lower parameters,
say, $b_s$, is of the form $q^{1+m}$, where $n$ and $m$ are
integers, with $n+m\ge 0$ if both conditions occur.

If the series does not terminate, it converges for $r<s$ and diverges for
$r>s$.  For $r=s$, it converges, provided it does not terminate from above,
for $|z|<1$ and, provided it does not terminate from below, for
$|b_1\cdots b_s/a_1\cdots a_rz|<1$.

\subsection{Summation and transformation formulae}

The most fundamental result in the theory of basic hypergeometric series
is the nonterminating $q$-binomial theorem
(cf.\ \cite[Equation (II.3)]{GR}),
\begin{equation}\label{1phi0}
{}_1\phi_0\!\left[\begin{matrix}
a\\-\end{matrix};q,z\right]=\frac{(az;q)_\infty}{(z;q)_\infty},
\qquad |z|<1.
\end{equation}

The following bilateral extension of \eqref{1phi0} is due to
Ramanujan (cf.\ \cite[Equation (II.29)]{GR}),
\begin{equation}\label{1psi1}
{}_1\psi_1\!\left[\begin{matrix}
a\\b\end{matrix};q,z\right]=
\frac{(q,az,q/az,b/a;q)_\infty}{(b,z,b/az,q/a;q)_\infty},
\qquad |b/a|<|z|<1.
\end{equation}
A particularly simple proof of \eqref{1psi1}
was given by Ismail~\cite{I}; it uses analyticity in the variable $b$ around
the origin, together with the identity theorem for analytic functions.  In this
paper we refer to this analytic-continuation argument as Ismail's argument.  Askey
and Ismail later used the same idea to give a short proof of Bailey's
very-well-poised ${}_6\psi_6$ summation~\cite{AI6}.  For an excellent survey
of Ramanujan's $_1\psi_1$ summation, see \cite{W}.

Among the many identities for basic hypergeometric series,
we recall several $_2\psi_2$ identities that will be used below.
In \cite[Eq.~(2.3)]{Ba}, Bailey derived the transformation
\begin{equation}\label{22tgl}
{}_2\psi_2\!\left[\begin{matrix}a,b\\
c,d\end{matrix};q,z\right]
=\frac{(az,d/a,c/b,dq/abz;q)_{\infty}}{(z,d,q/b,cd/abz;q)_{\infty}}\:
{}_2\psi_2\!\left[\begin{matrix}a,abz/d\\
az,c\end{matrix};q,\frac da\right],
\end{equation}
where $\max(|z|,|cd/abz|,|d/a|,|c/b|)<1$.  He iterated this transformation
to obtain \cite[Eq.~(2.4)]{Ba}
\begin{equation}\label{22tgl1}
{}_2\psi_2\!\left[\begin{matrix}a,b\\
c,d\end{matrix};q,z\right]
=\frac{(az,bz,cq/abz,dq/abz;q)_{\infty}}{(q/a,q/b,c,d;q)_{\infty}}\:
{}_2\psi_2\!\left[\begin{matrix}abz/c,abz/d\\
az,bz\end{matrix};q,\frac{cd}{abz}\right],
\end{equation}
where $\max(|z|,|cd/abz|)<1$.
Another important identity is Bailey's transformation of a general
$_2\psi_2$ series into a multiple of a very-well-poised $_6\psi_8$ series
(cf.\ \cite[Eq.~(3.2)]{Ba}; see also
\cite[Exercise~5.11, second identity]{GR}):
\begin{align}\label{bil68}\notag
{}_2\psi_2\!\left[\begin{matrix}
e,f\\aq/c,aq/d\end{matrix};q,\frac{aq}{ef}\right]
=\frac{(q/c,q/d,aq/e,aq/f;q)_\infty}
{(aq,q/a,aq/cd,aq/ef;q)_\infty}&\\*\times\;
{}_6\psi_8\!\left[\begin{matrix}qa^{\frac 12},-qa^{\frac 12},c,d,e,f\\
a^{\frac 12},-a^{\frac 12},aq/c,aq/d,aq/e,aq/f,0,0
\end{matrix};q,\frac{a^3q^2}{cdef}\right]&,
\end{align}
valid for $|aq/cd|<1$ and $|aq/ef|<1$.
We note that \eqref{22tgl} and \eqref{22tgl1} can be derived from
\eqref{bil68} by exploiting the symmetry of the $_6\psi_8$ series.

We will also use the following transformation, which writes a $_2\psi_2$
series as a sum of two multiples of $_2\phi_1$ series.  It was obtained by
Rosengren~\cite{Ro}; the form below is the equivalent formulation later given
by Chen, Chen, and Gu~\cite[Cor.~2.2]{CCG}:
\begin{align}\label{2psi2ccg}\notag
{}_2\psi_2\!\left[\begin{matrix}
a,b\\c,d\end{matrix};q,z\right]
=\frac{(q,c/b,q/d,abz/d,dq/abz;q)_\infty}
  {(q/a,q/b,c,az/d,cd/abz;q)_\infty}
{}_2\phi_1\!\left[\begin{matrix}
    cd/abz,d/a\\dq/az\end{matrix};q,\frac{bq}d\right]\\*
 + \frac{(q,az,q/az,b,d/a;q)_\infty}
  {(q/a,z,c,d,d/az;q)_\infty}
{}_2\phi_1\!\left[\begin{matrix}
    c/b,z\\azq/d\end{matrix};q,\frac{bq}d\right],
\end{align}
subject to $\max(|z|,|cd/abz|,|bq/d|)<1$.
As pointed out in \cite[Rem.~2.3]{CCG}, \eqref{2psi2ccg}
can also be formally obtained by taking a suitable limit in
\cite[Eq.~(III.34)]{GR} (which is a three-term $_3\phi_2$ identity
originally obtained by Sears in \cite[p.~175, Eq.~(10.2)]{S}).
The formal $c\to b$ limit of \eqref{2psi2ccg} gives (the $b\mapsto d$
case of) Ramanujan's ${}_1\psi_1$ summation \eqref{1psi1}.
We use \eqref{2psi2ccg} to establish the large $n$ and large $-n$
asymptotics of the bilateral $q$-ultraspherical functions introduced below.

For further material on basic hypergeometric series and, more generally, on
special functions, we refer to the textbooks by Gasper and Rahman~\cite{GR}
and by Andrews, Askey, and Roy~\cite{AAR}, respectively.  In our computations,
we repeatedly use elementary manipulations of $q$-shifted factorials; see
\cite[Appendix~I]{GR}.

\section{The continuous \texorpdfstring{$q$}{q}-ultraspherical polynomials}\label{sec:c}

\subsection{Definition and elementary properties}

We consider functions of $x=\cos\ta=(e^{\ri\ta}+e^{-\ri\ta})/2$
(where $\ta$ need not be real).

The continuous $q$-ultraspherical polynomials, which depend on
a parameter $\ba$ and the base $q$, are given by \cite[Sec.\ 13.2]{Ibook}
\begin{equation}\label{defC}
C_n(x;\ba\,|\,q)=\sum_{k=0}^n\frac{(\ba;q)_k\,(\ba;q)_{n-k}}
{(q;q)_k\,(q;q)_{n-k}}e^{\ri(n-2k)\ta}, \qquad x=\cos\ta.
\end{equation}
Rogers~\cite{R3} originally considered them in 1884, apparently unaware
of their orthogonality, in pursuit of what are now called the
Rogers--Ramanujan identities.
 
These functions, which can be written as
\begin{equation}\label{phyprep}
C_n(x;\ba\,|\,q)=\frac{(\ba;q)_n}{(q;q)_n}e^{\ri n\ta}\,
{}_2\phi_1\!\left[\begin{matrix}
\ba,q^{-n}\\q^{1-n}/\ba\end{matrix};q,qe^{-2\ri\ta}/\ba\right],
\end{equation}
are polynomials in $x$ of degree $n$.
They have the generating function
\begin{equation}
\sum_{n=0}^\infty C_n(x;\ba\,|\,q)t^n=
\frac{(\ba te^{\pm \ri\ta};q)_\infty}{(te^{\pm \ri\ta};q)_\infty},
\end{equation}
which readily follows from the definition \eqref{defC} and the
$q$-binomial theorem in \eqref{1phi0}.
From this generating function one readily deduces their
recurrence relation, which is \cite[Eq.~(13.2.12)]{Ibook}
\begin{align*}
&2x(1-\ba q^n)\,C_n(x;\ba\,|\,q)\\*
&=(1- q^{n+1})\,C_{n+1}(x;\ba\,|\,q)+
(1-\ba^2 q^{n-1})\,C_{n-1}(x;\ba\,|\,q),
\end{align*}
with the initial conditions
\begin{equation*}C_{-1}(x;\ba\,|\,q)=0\quad\;\text{and}\quad\;C_0(x;\ba\,|\,q)=1.
\end{equation*}
By induction, the $C_{2n+1}(x;\ba\,|\,q)$ are odd functions,
while the $C_{2n}(x;\ba\,|\,q)$ are even.
Their values at the origin are \cite[Eq.~(13.2.19)]{Ibook}
\begin{equation*}
C_{2n+1}(0;\ba\,|\,q)=0\quad\,\text{and}\quad\,
C_{2n}(0;\ba\,|\,q)=(-1)^n\frac{(\ba^2;q^2)_n}{(q^2;q^2)_n},
\end{equation*}
for all nonnegative integers $n$.

Let $\mathcal D_q$ denote the \textit{Askey--Wilson operator},
defined on functions of $x=(z+1/z)/2$
(recall that $x=\cos\theta$, so $z=e^{\ri\ta}$) by \cite[Eq.~(12.1.10)]{Ibook}
\begin{align}
\mathcal D_q f(x)&=\mathcal D_q f((z+1/z)/2)\notag\\
:\!\!\,\!&=\frac{f((q^{\frac 12}z+q^{-\frac 12}/z)/2)-
f((q^{-\frac 12}z+q^{\frac 12}/z)/2)}
{((q^{\frac 12}z+q^{-\frac 12}/z)/2)-((q^{-\frac 12}z+q^{\frac 12}/z)/2)}\notag\\[.3em]
&=\frac{f((q^{\frac 12}z+q^{-\frac 12}/z)/2)-
f((q^{-\frac 12}z+q^{\frac 12}/z)/2)}
{(q^{\frac 12}-q^{-\frac 12})(z-1/z)/2}.
\end{align}
The action of $\mathcal D_q$ on the continuous $q$-ultraspherical polynomials
is \cite[Eq.~(13.2.23)]{Ibook}
\begin{equation*}
\mathcal D_q C_n(x;\ba\,|\,q)=
\frac{2(1-\ba)}{(1-q)}q^{\frac{1-n}2}
C_{n-1}(x;q\ba\,|\,q).
\end{equation*}

\subsection{Orthogonality and integrals}

As was established by Askey and Ismail \cite{AI} (see also
\cite[Thm.~13.2.1]{Ibook}), the continuous $q$-ultraspherical polynomials
satisfy, for $|\ba|<1$, the orthogonality relation
\begin{subequations}
\label{orth}
\begin{align}\label{orthrel}
&\frac 1{2\pi}\int_{-1}^1C_m(x;\ba\,|\,q)\,C_n(x;\ba\,|\,q)\,
w(x\,|\,\ba){\rd x}\notag\\
&\quad=
\frac{(\ba,q\ba;q)_\infty}{(q,\ba^2;q)_\infty}
\frac{(\ba^2;q)_n}{(q;q)_n}\frac{(1-\ba)}{(1-\ba q^n)}\,
\da_{m,n},
\end{align}
where
\begin{equation}\label{orthw}
 w(x\,|\,\ba)=\frac{(e^{2\ri\ta},e^{-2\ri\ta};q)_\infty}
  {(\ba e^{2\ri\ta},\ba e^{-2\ri\ta};q)_\infty}
\frac 1{\sqrt{1-x^2}},\quad\;\text{$x=\cos\ta$, $0\le \ta\le \pi$,}
\end{equation}
\end{subequations}
is the corresponding weight function.

Askey and Ismail \cite[p.~6]{AI} noted that the corresponding
positive definiteness condition for the full orthogonality functional is
equivalent to
$$
\frac{(1-\beta^2 q^{n-1})(1-q^n)}{(1-\beta q^{n-1})(1-\beta q^n)}>0,\qquad
\text{for $n=1,2,\ldots$,}
$$
from which one determines the following restrictions on
the real parameters $q$ and $\beta$:
\begin{alignat*}{2}
-1<{}&\ba<q^{-\frac 12}&\qquad\qquad&\text{for $0<q<1$,}\\*
-1<{}&\ba<-q^{-1}&\qquad\qquad&\text{for $-1<q<0$},\\*
&\ba>-1&\quad&\text{for $q=0$},\\*
&\ba=q^\lambda\quad\text{with $\lambda>-\frac 12$}&\qquad&\text{for $q\to 1$},\\*
&\ba=|q|^\lambda\quad\text{with $\lambda>-1$}&\qquad\qquad& \\*
or\quad &\ba=-|q|^\lambda\quad\text{with $\lambda>0$}&\qquad&\text{for $q\to -1$}.
\end{alignat*}
For later comparison, we also record the corresponding finite mass-point
form of the continuous $q$-ultraspherical orthogonality relation, obtained
from the Askey--Wilson orthogonality relation when poles have crossed the
unit circle (cf.\ \cite[Ch.~14]{KLS}).  In particular, if $0<q<1$ and
$\ba>1$, set
\begin{subequations}\label{masspoints}
\begin{align}
x_j&=\frac{(\ba q^j)^{1/2}+(\ba q^j)^{-1/2}}2,
\qquad j=0,1,2,\ldots,\label{masspoints-x}\\*
\lambda_j&=
\frac12\,
\frac{(\ba^{-1},q \ba;q)_\infty}{(q,\ba^2;q)_\infty}
\frac{(1-\ba q^j)(\ba^2;q)_j}{(q;q)_j}\,\ba^{-j}.
\label{masspoints-lambda}
\end{align}
\end{subequations}
Then, away from exceptional colliding-pole values (where limits are to be
understood),
\begin{align}\label{orthmass}
&\frac 1{2\pi}\int_{-1}^1
C_m(x;\ba\,|\,q)C_n(x;\ba\,|\,q)w(x\,|\,\ba)\rd x\notag\\*
&\quad+
\sum_{\substack{j\ge0\\ \ba q^j>1}}\lambda_j
\Big(C_m(x_j;\ba\,|\,q)C_n(x_j;\ba\,|\,q)
+C_m(-x_j;\ba\,|\,q)C_n(-x_j;\ba\,|\,q)\Big)\notag\\
&=\frac{(\ba,q\ba;q)_\infty}{(q,\ba^2;q)_\infty}
\frac{(\ba^2;q)_n}{(q;q)_n}\frac{(1-\ba)}{(1-\ba q^n)}\,
\da_{m,n}.
\end{align}
Ismail~\cite[Thm.~13.2.2]{Ibook} further observed that the orthogonality
relation in \eqref{orth} is equivalent to the following integral evaluation.
\begin{equation}\label{thm1322}
\frac 1{2\pi}\int_{-1}^1
\frac{(\ba t_1e^{\pm\ri\ta},\ba t_2e^{\pm\ri\ta},e^{\pm 2\ri\ta};q)_\infty}
{(t_1e^{\pm\ri\ta},t_2e^{\pm\ri\ta},\ba e^{\pm 2\ri\ta};q)_\infty}
\frac{\rd x}{\sqrt{1-x^2}}
=\frac{(\ba,q\ba;q)_\infty}{(q,\ba^2;q)_\infty}\,
{}_2\phi_1\!\left[\begin{matrix}
\ba^2,\ba\\q\ba\end{matrix};q,t_1t_2\right],
\end{equation}
where $|t_1|<1$ and $|t_2|<1$.
In fact, this result is an immediate consequence of combining \eqref{orth}
with the well-known generating function for the continuous $q$-ultraspherical
polynomials~\cite[Eq.~(13.2.8)]{Ibook}
\begin{equation}\label{gfuni}
\sum_{n=0}^\infty C_n(x;\ba\,|\,q)t^n=
\frac{(\ba te^{\pm\ri\ta};q)_\infty}
{(te^{\pm\ri\ta};q)_\infty},
\end{equation}
where $|te^{\pm\ri\ta}|<1$.  Multiplying both sides of \eqref{orthrel}
by $t_1^mt_2^n$ and summing over all $m,n\ge 0$ using \eqref{gfuni},
one readily obtains \eqref{thm1322}.
Conversely, taking coefficients of $t_1^mt_2^n$ on both sides of \eqref{thm1322}
gives \eqref{orthrel}. If one uses the second generating function
for the continuous $q$-ultraspherical polynomials in Corollary~\ref{corgf2},
then, under the explicit pure-integral hypotheses
\begin{equation*}
|\ba|<1,\qquad |t_1|<1,\qquad |t_2|<1,
\end{equation*}
and away from zeros of the displayed denominator factors, one obtains the
following two pure-integral contiguous companions to \eqref{thm1322}.  Their
right-hand sides factorize completely into products of linear factors:
\begin{subequations}
\begin{align}
&\frac 1{2\pi}\int_{-1}^1
\frac{(\ba t_1e^{\pm\ri\ta},
q\ba t_2e^{\pm\ri\ta},
e^{\pm 2\ri\ta};q)_\infty}
{(t_1e^{\pm\ri\ta},
t_2e^{\pm\ri\ta},
\ba e^{\pm 2\ri\ta};q)_\infty}
\frac{\rd x}{\sqrt{1-x^2}}\notag\\*
&=\frac{(\ba,q\ba,
\ba^2t_1t_2;q)_\infty}
  {(q,\ba^2,
 t_1t_2;q)_\infty}
 \frac 1{(1-\ba t_2^2)},\label{orth-idgf2-cont}
\end{align}
\begin{align}
&\frac 1{2\pi}\int_{-1}^1
\frac{(q\ba t_1e^{\pm\ri\ta},
q\ba t_2e^{\pm\ri\ta},
e^{\pm 2\ri\ta};q)_\infty}
{(t_1e^{\pm\ri\ta},
t_2e^{\pm\ri\ta},
\ba e^{\pm 2\ri\ta};q)_\infty}
\frac{\rd x}{\sqrt{1-x^2}}\notag\\*
&=\frac{(\ba,q\ba,
q\ba^2t_1t_2;q)_\infty}
  {(q,\ba^2,
 t_1t_2;q)_\infty}
 \frac {(1+\ba t_1t_2)}{(1-\ba t_1^2)(1-\ba t_2^2)}.\label{orth-idgf3-cont}
\end{align}
\end{subequations}
Indeed, to prove \eqref{orth-idgf2-cont}, multiply both sides of \eqref{orthrel}
by $(1-\ba q^n)t_1^mt_2^n$ and sum over all $m,n\ge 0$.
The sums can be evaluated by \eqref{gfuni} and Corollary~\ref{corgf2}.
Division of both sides of the identity by $(1-\ba)(1-\ba t_2^2)$
and simplification of the identity using the $(a;z)\mapsto (\ba^2,t_1t_2)$
instance of the nonterminating $q$-binomial theorem \eqref{1phi0}
gives \eqref{orth-idgf2-cont}.
To prove \eqref{orth-idgf3-cont}, multiply both sides of \eqref{orthrel}
by $(1-\ba q^m)(1-\ba q^n)t_1^mt_2^n$ and sum over all $m,n\ge 0$.
The sums can be evaluated by Corollary~\ref{corgf2}.
Division of both sides of the identity by
$(1-\ba)^2(1-\ba t_1^2)(1-\ba t_2^2)$
and simplification of the identity using a weighted sum of
two instances of \eqref{1phi0}
(namely $(a;z)\mapsto (\ba^2,t_1t_2)$ and $(a;z)\mapsto (\ba^2,qt_1t_2)$)
gives \eqref{orth-idgf3-cont}.

\begin{remark}
The displayed pure-integral evaluations should not be read as unrestricted
orthogonality consequences.  Outside the mass-free range the full continuous
orthogonality relation is \eqref{orthmass}, and the mass terms produced by the
generating-function summation must be included.  The corresponding bilateral
integral evaluations with mass aggregates are stated in
Theorem~\ref{thm:mass-integral-evals}.
\end{remark}


\subsection{Rogers' linearization formula}

One of the remarkable properties of the continuous
$q$-ultraspherical polynomials is that they possess a linearization formula
whose coefficients completely factorize.
This result is due to Rogers~\cite{R3} (cf.\ also \cite[Thm.~13.3.2]{Ibook}),
who obtained the following formula:
\begin{align}\label{linear}
&C_m(x;\ba\,|\,q)\,C_n(x;\ba\,|\,q)\notag\\
&=\sum_{k=0}^{\min(m,n)}
\bigg(\frac{(q;q)_{m+n-2k}(\ba;q)_{m-k}(\ba;q)_{n-k}(\ba;q)_k(\ba^2;q)_{m+n-k}}
{(\ba^2;q)_{m+n-2k}(q;q)_{m-k}(q;q)_{n-k}(q;q)_k(q\ba;q)_{m+n-k}}\notag\\*
&\qquad\qquad\quad\times\frac{(1-\ba q^{m+n-2k})}{(1-\ba)}\,C_{m+n-2k}(x;\ba\,|\,q)\bigg).
\end{align}
A mixed bilateral extension of Rogers' linearization formula \eqref{linear},
together with a quasi-linearized companion, will be given in
Section~\ref{sec:linear}.

\section{The bilateral \texorpdfstring{$q$}{q}-ultraspherical functions}\label{sec:bilf}

\subsection{Definition and bilateral series representations}

For two parameters $\ba,\ga$, base $q$, and variable $x=\cos\ta$,
we define the {\em bilateral $q$-ultraspherical functions} by
\begin{equation}\label{defbilC}
C_n(x;\ba,\ga\,|\,q):=
\sum_{k=-\infty}^\infty\frac{(\ba;q)_k(\ba;q)_{n-k}}
{(q\ga;q)_k(q\ga;q)_{n-k}}e^{\ri(n-2k)\ta},
\end{equation}
for $|q\ga e^{\pm 2\ri\ta}/\ba|<1$.
(The parameters can be further restricted for a positive measure when needed.)
It is clear that, as $\ga\to 1$, the bilateral $q$-ultraspherical function
$C_n(x;\ba,\ga\,|\,q)$ reduces to the continuous
$q$-ultraspherical polynomial $C_n(x;\ba\,|\,q)$ in \eqref{defC}.
Also, shifting the summation index $k\mapsto k+n$ shows that the symmetry
\begin{equation}\label{sym}
C_{n}(x;\ba,\ga\,|\,q)=\bigg(\frac{\ba}{q\ga}\bigg)^nC_{-n}(x;1/\ga,1/\ba\,|\,q)
\end{equation}
holds for all $n\in\mathbb Z$.

The function $C_n(x;\ba,\ga\,|\,q)$ possesses the following
representations as multiples of a convergent
bilateral basic hypergeometric series:
\begin{subequations}\label{hypreps}
\begin{align}
&C_n(x;\ba,\ga\,|\,q)\notag\\
&=\frac{(\ba;q)_n}{(q\ga;q)_n}e^{\ri n\ta}\,
{}_2\psi_2\!\left[\begin{matrix}\ba,q^{-n}/\ga\\
    q\ga,q^{1-n}/\ba\end{matrix};q,q\ga e^{-2\ri\ta}/\ba\right]\label{hyprepa}\\
 &=\frac{(\ba;q)_n}{(q\ga e^{2\ri\ta};q)_n}e^{\ri n\ta}
   \frac{(q\ga/\ba,q\ga/\ba,q\ga e^{\pm 2\ri\ta};q)_\infty}
   {(q\ga,q\ga,q\ga e^{\pm 2\ri\ta}/\ba;q)_\infty}
   {}_2\psi_2\!\left[\begin{matrix}\ba,q^{-n}e^{-2\ri\ta}/\ga\\
       q\ga e^{-2\ri\ta},q^{1-n}/\ba\end{matrix}
  ;q,q\ga/\ba\right]
  \label{hyprepb}\\
 &=\frac{(\ba e^{2\ri\ta};q)_n}{(q\ga;q)_n}e^{-\ri n\ta}
   \frac{(q\ga/\ba,q\ga/\ba,q e^{\pm 2\ri\ta}/\ba;q)_\infty}
   {(q/\ba,q/\ba,q\ga e^{\pm 2\ri\ta}/\ba;q)_\infty}
   {}_2\psi_2\!\left[\begin{matrix}\ba e^{-2\ri\ta},q^{-n}/\ga\\
      q\ga,q^{1-n}e^{-2\ri\ta}/\ba\end{matrix}
  ;q,q\ga/\ba\right]
  \label{hyprepc}\\
 &=\frac{(\ba e^{2\ri\ta};q)_n}{(q\ga e^{2\ri\ta};q)_n}e^{-\ri n\ta}
   \frac{(q\ga e^{\pm 2\ri\ta},
   q e^{\pm 2\ri\ta}/\ba;q)_\infty}
   {(q/\ba,q/\ba,q\ga,q\ga;q)_\infty}
   {}_2\psi_2\!\left[\begin{matrix}\ba e^{-2\ri\ta},q^{-n}e^{-2\ri\ta}/\ga\\
       q\ga e^{-2\ri\ta},q^{1-n}e^{-2\ri\ta}/\ba\end{matrix}
  ;q,q\ga e^{2\ri\ta}/\ba\right]
  \label{hyprepd}\\
  &=\frac{(\ba,\ba e^{2\ri\ta};q)_n}
    {(q\ga,q\ga e^{2\ri\ta};q)_n}
    \frac{q^ne^{\ri (n+1)\ta}}{2\ri\sin\ta\,\ba^n}
    \frac{(q\ga e^{\pm 2\ri\ta},
    qe^{\pm 2\ri\ta}/\ba;q)_\infty}
    {(qe^{\pm 2\ri\ta},
    q\ga e^{\pm 2\ri\ta}/\ba;q)_\infty}\notag\\
&\quad\;\times\sum_{k=-\infty}^\infty
  \frac{(\ba,\ba e^{-2\ri\ta},q^{-n}/\ga,q^{-n}e^{-2\ri\ta}/\ga;q)_k}
 {(q^{1-n}e^{-2\ri\ta}/\ba,q^{1-n}/\ba,q\ga e^{-2\ri\ta},q\ga;q)_k}
 \big(1-q^{-n+2k}e^{-2\ri\ta}\big)
 \bigg(\frac{\ga^2e^{-2\ri\ta}}{q^{n-1}\ba^2}\bigg)^kq^{k^2}
 .\label{hyprepe}
\end{align}
\end{subequations}
Since the function $C_n(x;\ba,\ga\,|\,q)$ is symmetric in $e^{\ri\ta}$
and $e^{-\ri\ta}$, we may replace  $e^{\ri\ta}$ by $e^{-\ri\ta}$ in
any of the bilateral series representations in \eqref{hypreps} to get
additional representations for $C_n(x;\ba,\ga\,|\,q)$.

The $_2\psi_2$ series appearing in \eqref{hyprepa}--\eqref{hyprepd}
are all \textit{well-poised}.
(See \cite[p.~39 and p.~138]{GR} for this terminology in the setting of
basic hypergeometric and bilateral basic hypergeometric series.)
The first series representation
in \eqref{hyprepa} directly stems from the defining relation
\eqref{defbilC} for $C_n(x;\ba,\ga\,|\,q)$.
The representations in \eqref{hyprepb} and \eqref{hyprepc} are
obtained by applying \eqref{22tgl} to \eqref{hyprepa}.
The representation in \eqref{hyprepd} is obtained by
applying \eqref{22tgl1} to \eqref{hyprepa}.
Finally, application of \eqref{bil68} to \eqref{hyprepa}
(or to any of the other representations of $C_n(x;\ba,\ga\,|\,q)$
as a multiple of a $_2\psi_2$ series) yields the bilateral series
representation in \eqref{hyprepe}.

In the two ${}_2\psi_2$ series in \eqref{hyprepa} and \eqref{hyprepd}
we require $|q\ga e^{\pm 2\ri\ta}/\ba|<1$,
while in \eqref{hyprepb} and \eqref{hyprepc}
we only require $|q\ga/\ba|<1$ (which is a larger region than
$|q\ga e^{\pm 2\ri\ta}/\ba|<1$ if $\ta$ is not a real number),
for absolute convergence.
The bilateral series in \eqref{hyprepe} converges absolutely
everywhere (because of the quadratic powers of $q$ appearing
as a factor in the summand; throughout we assume $|q|<1$)
and imposes no restrictions on the parameters, apart from avoiding poles.
Thus, the bilateral series in \eqref{hyprepe} gives a formula for
\textit{analytic continuation} of the bilateral $q$-ultraspherical
function $C_n(x;\ba,\ga\,|\,q)$ beyond the region
$|q\ga e^{\pm 2\ri\ta}/\ba|<1$,
to the full complex plane avoiding the set of poles.
In addition, the series representation in \eqref{hyprepe} is practical
for computational purposes, as the series converges quickly.

We point out that applications of Bailey's transformation
\eqref{22tgl} to \eqref{hyprepa} that result in the argument of
the ${}_2\psi_2$ series containing the factor $q^n$ or $q^{-n}$
are less suitable for the present list, since they do not provide a
single absolutely convergent representation valid uniformly for all
$n\in\mathbb Z$.
For fixed $n$, or when one studies a one-sided limit such as
$n\to\infty$ or $n\to-\infty$, such transformed representations may
nevertheless be useful,
for instance in asymptotic considerations or as meromorphic continuation
formulae between overlapping regions of convergence.

We also have the following useful representations of $C_n(x;\ba,\ga\,|\,q)$
as a sum of two multiples of ${}_2\phi_1$ series with manifest
$e^{\pm\ri\ta}$ symmetry.
\begin{subequations}\label{2term}
\begin{align}
&C_n(x;\ba,\ga\,|\,q)\notag\\
&=\frac{(q,\ba,q\ga/\ba;q)_\infty}
{(q\ga,q\ga,q/\ba;q)_\infty}\notag\\*
&\quad\;\times\left(\frac{(\ba e^{-2\ri\ta},q e^{2\ri\ta}/\ba;q)_\infty}
{(e^{-2\ri\ta},q\ga e^{2\ri\ta}/\ba;q)_\infty}e^{\ri n\ta}
{}_2\phi_1\!\left[\begin{matrix}q\ga/\ba,q\ga e^{2\ri\ta}/\ba\\
q e^{2\ri\ta}\end{matrix};q,\ba^2 q^n\right]
\right.\notag\\*
&\qquad\quad +\left.
\frac{(\ba e^{2\ri\ta},q e^{-2\ri\ta}/\ba;q)_\infty}
{(e^{2\ri\ta},q\ga e^{-2\ri\ta}/\ba;q)_\infty}e^{-\ri n\ta}
{}_2\phi_1\!\left[\begin{matrix}q\ga/\ba,q\ga e^{-2\ri\ta}/\ba\\
q e^{-2\ri\ta}\end{matrix};q,\ba^2q^n\right]
\right)\label{2terma}\\
&=\frac{(q,\ba,q\ga/\ba,q\ga^2;q)_\infty}
{(q\ga,q\ga,q/\ba,\ba^2;q)_\infty}
\frac{(\ba^2;q)_n}{(q\ga^2;q)_n}\notag\\*
&\quad\;\times\left(\frac{(\ba e^{-2\ri\ta},q e^{2\ri\ta}/\ba;q)_\infty}
{(e^{-2\ri\ta},q\ga e^{2\ri\ta}/\ba;q)_\infty}e^{\ri n\ta}
{}_2\phi_1\!\left[\begin{matrix}\ba/\ga,\ba e^{2\ri\ta}/\ga\\
q e^{2\ri\ta}\end{matrix};q,\ga^2q^{n+1}\right]
\right.\notag\\*
&\qquad\quad +\left.
\frac{(\ba e^{2\ri\ta},q e^{-2\ri\ta}/\ba;q)_\infty}
{(e^{2\ri\ta},q\ga e^{-2\ri\ta}/\ba;q)_\infty}e^{-\ri n\ta}
{}_2\phi_1\!\left[\begin{matrix}\ba/\ga,\ba e^{-2\ri\ta}/\ga\\
q e^{-2\ri\ta}\end{matrix};q,\ga^2q^{n+1}\right]
\right)\label{2termb}
\end{align}
\end{subequations}

The expression in \eqref{2terma} is a direct consequence of
\eqref{2psi2ccg} applied to \eqref{hyprepa}.
The further application of the $q$-Euler transformation~\cite[Eq.~(III.3)]{GR}
to the two ${}_2\phi_1$ series in \eqref{2terma} gives \eqref{2termb}.
The latter formula extends~\cite[Eq.~(1.13)]{RV},
to which \eqref{2termb} reduces for $\ga\to 1$.
Since $q^n$ appears only in the argument
but not in the parameters of the above ${}_2\phi_1$ series, the
series representations in \eqref{2term} are useful
for finding the large $n$ (and, using
\eqref{sym}, the large $-n$)
asymptotics for the bilateral $q$-ultraspherical
functions $C_n(x;\ba,\ga\,|\,q)$.  We provide full
details on their asymptotics in Section~\ref{sec:asym}.

\begin{remark}
There are several possible ways to extend polynomials
to bilateral series.
One may ask what our motivation is to bilaterally extend the continuous
$q$-ultraspherical polynomials precisely as in \eqref{defbilC}.
The answer lies in the useful properties that carry over
from the unilateral to the bilateral case.  In particular,
the results in Theorems~\ref{thmgf} and \ref{thmrec} are especially simple.
Furthermore, it is interesting to see that the recurrence relation \eqref{rec}
satisfied by the bilateral $q$-ultraspherical
functions $C_n(x;\ba,\ga\,|\,q)$ has the same coefficient form as the
recurrence relation for the \textit{associated}
continuous $q$-ultraspherical polynomials $C_n^{(a)}(x;b\,|\,q)$
(see \cite[Eq.\ (2.2)]{BI}), after the substitution
$a=\ga^2$ and $b=\ba/\ga$.  The initial conditions are different, however.
This is therefore a comparison at the level of the
recurrence only; it does not identify the bilateral functions with the
associated polynomials, nor does it transfer their orthogonality.  The
orthogonality functionals constructed in Section~\ref{sec:orth} are obtained
from the bilateral generating functions and terminating specializations.
This connection suggests that it should make sense to study bilateral
extensions of other families of orthogonal polynomials that have associated
companions,
in particular in such a way that the
recurrence relation for the respective associated family
is also the underlying recurrence relation for the bilateral extension.
We leave this for future research.
\end{remark}

\subsection{Generating functions, recurrence, and special values}

We have the following bilateral generating function for the
bilateral $q$-ultraspherical functions:
\begin{theorem}\label{thmgf}
We have
\begin{equation}\label{gfbil}
\sum_{n=-\infty}^\infty C_n(x;\ba,\ga\,|\,q)t^n=
\frac{(q,q\ga/\ba;q)_\infty^2}{(q\ga,q/\ba;q)_\infty^2}
\frac{(\ba te^{\pm\ri\ta},qe^{\pm\ri\ta}/\ba t;q)_\infty}
{(te^{\pm\ri\ta},q\ga e^{\pm\ri\ta}/\ba t;q)_\infty},
\end{equation}
where $|q\ga/\ba|<|te^{\pm\ri\ta}|<1$.
\end{theorem}
\begin{proof}
We have
\begin{align*}
\sum_{n=-\infty}^\infty C_n(x;\ba,\ga\,|\,q)t^n&=
\sum_{n,\,k=-\infty}^\infty\frac{(\ba;q)_k}{(q\ga;q)_k}
\frac{(\ba;q)_{n-k}}{(q\ga;q)_{n-k}}e^{\ri(-k+(n-k))\ta}t^{k+(n-k)}\\
&=\sum_{k=-\infty}^\infty\frac{(\ba;q)_k}{(q\ga;q)_k}e^{-\ri k\ta}t^k\,
\sum_{n=-\infty}^\infty\frac{(\ba;q)_{n}}{(q\ga;q)_{n}}e^{\ri n\ta}t^{n}\\
&=\frac{(q,\ba te^{-\ri\ta},qe^{\ri\ta}/\ba t,q\ga/\ba;q)_\infty}
{(q\ga,te^{-\ri\ta},q\ga e^{\ri\ta}/\ba t,q/\ba;q)_\infty}
\frac{(q,\ba te^{\ri\ta},qe^{-\ri\ta}/\ba t,q\ga/\ba;q)_\infty}
{(q\ga,te^{\ri\ta},q\ga e^{-\ri\ta}/\ba t,q/\ba;q)_\infty}.
\end{align*}
Here the last equality follows from two applications of Ramanujan's
$_1\psi_1$ summation \eqref{1psi1}.
\end{proof}

The $t\mapsto qt$ case of \eqref{gfbil} is
\begin{subequations}\label{gfbil-qt}
\begin{align}
\sum_{n=-\infty}^{\infty} C_n(x;\ba,\ga\,|\,q)q^nt^n&=
\frac{(q,q\ga/\ba;q)_\infty^2}{(q\ga,q/\ba;q)_\infty^2}
\frac{(q\ba te^{\pm\ri\ta},e^{\pm\ri\ta}/\ba t;q)_\infty}
{(qt e^{\pm\ri\ta},\ga e^{\pm\ri\ta}/\ba t;q)_\infty}\notag\\
&=\frac{(1-te^{\ri\ta})(1-te^{-\ri\ta})
(1-e^{\ri\ta}/\ba t)(1-e^{-\ri\ta}/\ba t)}
{(1-\ba te^{\ri\ta})(1-\ba te^{-\ri\ta})
(1-\ga e^{\ri\ta}/\ba t)(1-\ga e^{-\ri\ta}/\ba t)}\notag\\*
&\quad\;\times
\frac{(q,q\ga/\ba;q)_\infty^2}{(q\ga,q/\ba;q)_\infty^2}
\frac{(\ba te^{\pm\ri\ta},qe^{\pm\ri\ta}/\ba t;q)_\infty}
{(te^{\pm\ri\ta},q\ga e^{\pm\ri\ta}/\ba t;q)_\infty},
\end{align}
where $t$ is in the annulus $|\ga e^{\pm\ri\ta}/\ba|<|t|<|e^{\pm\ri\ta}/q|$.
Now since
\begin{equation}
\frac{(1-e^{\ri\ta}/\ba t)(1-e^{-\ri\ta}/\ba t)}
{(1-\ba te^{\ri\ta})(1-\ba te^{-\ri\ta})}=\frac 1{\ba^2t^2},
\end{equation}
\end{subequations}
the right-hand side simplifies and (recalling $x=(e^{\ri\ta}+e^{-\ri\ta})/2$)
we deduce
\begin{align}
&\sum_{n=-\infty}^{\infty} C_n(x;\ba,\ga\,|\,q)q^nt^n\notag\\
&=
\frac{(1-te^{\ri\ta})(1-te^{-\ri\ta})}
{\ba^2t^2(1-\ga e^{\ri\ta}/\ba t)(1-\ga e^{-\ri\ta}/\ba t)}\,
\frac{(q,q\ga/\ba;q)_\infty^2}{(q\ga,q/\ba;q)_\infty^2}
\frac{(\ba te^{\pm\ri\ta},qe^{\pm\ri\ta}/\ba t;q)_\infty}
{(te^{\pm\ri\ta},q\ga e^{\pm\ri\ta}/\ba t;q)_\infty}\notag\\
&=
\frac{1-2tx+t^2}{\ga^2-2\ba\ga xt+\ba^2t^2}\,
\frac{(q,q\ga/\ba;q)_\infty^2}{(q\ga,q/\ba;q)_\infty^2}
\frac{(\ba te^{\pm\ri\ta},qe^{\pm\ri\ta}/\ba t;q)_\infty}
{(te^{\pm\ri\ta},q\ga e^{\pm\ri\ta}/\ba t;q)_\infty}\notag\\
&=\label{recrel}
\frac{1-2tx+t^2}{\ga^2-2\ba\ga xt+\ba^2t^2}\,
\sum_{n=-\infty}^{\infty} C_n(x;\ba,\ga\,|\,q)t^n,
\end{align}
where we now assume that $t$ lies in the intersection of the annuli
$|\ga e^{\pm\ri\ta}/\ba|<|t|<|e^{\pm\ri\ta}/q|$ and
$|q\ga e^{\pm\ri\ta}/\ba|<|t|<|e^{\pm\ri\ta}|$, i.e., where
$|\ga e^{\pm\ri\ta}/\ba|<|t|<|e^{\pm\ri\ta}|$.
Multiplying the left- and right-hand sides of \eqref{recrel} by
$\ga^2-2\ba\ga xt+\ba^2t^2$ and comparing coefficients of $t^n$,
we arrive at the following result (where we recall from \eqref{defbilC}
the condition $|q\ga e^{\pm 2\ri\ta}/\ba|<1$ that is required for
convergence of the bilateral $q$-ultraspherical functions;
this condition is consistent with the above annuli being non-empty):
\begin{theorem}\label{thmrec}
The bilateral $q$-ultraspherical functions satisfy for all $n\in\Z$
the following recurrence relation:
\begin{align}
&2x(1-\ba\ga q^{n})\,C_n(x;\ba,\ga\,|\,q)\notag\\*
&=(1-\ga^2 q^{n+1})\,C_{n+1}(x;\ba,\ga\,|\,q)+
(1-\ba^2q^{n-1})\,C_{n-1}(x;\ba,\ga\,|\,q).\label{rec}
\end{align}
Moreover, they are uniquely determined by \eqref{rec} and
the two initial conditions
\begin{subequations}
\begin{align}
\label{C0spec}
C_0(x;\ba,\ga\,|\,q)&={}_2\psi_2\!\left[\begin{matrix}
\ba,1/\ga\\ q\ga,q/\ba\end{matrix};q,q\ga e^{-2\ri\ta}/\ba\right],\\*
C_{-1}(x;\ba,\ga\,|\,q)&=\frac{(1-\ga)}{(1-\ba/q)}e^{-\ri\ta}\,
{}_2\psi_2\!\left[\begin{matrix}
\ba,q/\ga\\ q\ga,q^2/\ba\end{matrix};q,q\ga e^{-2\ri\ta}/\ba\right].
\end{align}
\end{subequations}
\end{theorem}

\begin{remark}
In the classical unilateral case, it is principally clear from
the basic hypergeometric representation in \eqref{phyprep},
together with standard contiguous relations~\cite{Kr},
that the continuous $q$-ultraspherical polynomials
satisfy a three-term relation.  However, it is not a priori clear
that the three-term relation should have the same form
as the one compatible with Favard's theorem~\cite[Sec.~3.2]{Sz},
with $x$ appearing only in the coefficient of the term $C_n$.
In the bilateral case, the existence of the
three-term relation does not straightforwardly follow from the fact that
the bilateral $q$-ultraspherical functions
are multiples of $_2\psi_2$ series.
The reader should be reminded that contiguous relations
implicitly make use of identities such as
$(1-uq^k)/(1-u)=q^k+(1-q^k)/(1-u)$ which are connected
to splitting a unilateral basic hypergeometric sum in two parts and
incrementing the summation index $k$ in the second sum (since the
factor $(1-q^k)$ kills the $k=0$ term).
The existence of the three-term recurrence relation \eqref{rec}
is due to the very specific form of the basic
hypergeometric representation in \eqref{hyprepa}.
\end{remark}

The bilateral $q$-ultraspherical functions $C_0(x)$ and $C_{-1}(x)$
in general do not evaluate in closed form; however, from the $b=\ba$,
$c=1/\ga$, $d=q^{\frac 12}$ case of
\begin{subequations}
\begin{equation}
{}_3\psi_3\!\left[\begin{matrix}
b,c,d\\q/b,q/c,q/d\end{matrix};q,\frac q{bcd}\right]=
\frac{(q,q/bc,q/bd,q/cd;q)_\infty}{(q/b,q/c,q/d,q/bcd;q)_\infty},
\end{equation}
and the $b=\ba$, $c=q/\ga$, $d=-q$ case
of
\begin{equation}
{}_3\psi_3\!\left[\begin{matrix}
b,c,d\\q^2/b,q^2/c,q^2/d\end{matrix};q,\frac{q^2}{bcd}\right]=
\frac{(q,q^2/bc,q^2/bd,q^2/cd;q)_\infty}
{(q^2/b,q^2/c,q^2/d,q^2/bcd;q)_\infty},
\end{equation}
\end{subequations}
which are formulae by Bailey for specific well-poised
$_3\psi_3$ series (cf.\ \cite[Ex.~5.18~(i) and (ii)]{GR}),
we can deduce the following evaluations for special $x$:
\begin{subequations}
\begin{align}
  C_0\big((q^{\frac 14}+q^{-\frac 14})/2;\ba,\ga\,|\,q\big)
  &=
\frac{(q,q\ga/\ba,q^{\frac 12}\ga,q^{\frac 12}/\ba;q)_\infty}
{(q^{\frac 12},q^{\frac 12}\ga/\ba,q\ga,q/\ba;q)_\infty},\\*
  C_{-1}\big(\ri(q^{\frac 12}-q^{-\frac 12})/2;\ba,\ga\,|\,q\big)
  &=\frac{\ri\, q^{\frac 12}(1-\ga)}{\ba}
 \frac{(q,q\ga/\ba,-\ga,-q/\ba;q)_\infty}
{(-q,-\ga/\ba,q\ga,q/\ba;q)_\infty}.
\end{align}
\end{subequations}

The $C_{2n+1}(x;\ba,\ga\,|\,q)$ are odd functions,
the $C_{2n}(x;\ba,\ga\,|\,q)$ are even.

Their constant terms are
\begin{align*}
C_{2n+1}(0;\ba,\ga\,|\,q)&=0,\\
\intertext{and}
C_{2n}(0;\ba,\ga\,|\,q)&=(-1)^n\frac{(\ba^2;q^2)_n}{(q^2\ga^2;q^2)_n}\,
\frac{(q,q\ga/\ba,-q\ga,-q/\ba;q)_\infty}{(-q,-q\ga/\ba,q\ga,q/\ba;q)_\infty},
\end{align*}
for all $n\in\mathbb Z$, by virtue of an instance of the bilateral $q$-Kummer
summation \cite[Eq.~(II.30)]{GR}, namely
\begin{equation}
{}_2\psi_2\!\left[\begin{matrix}
b,c\\aq/b,aq/c\end{matrix};q,-\frac{aq}{bc}\right]=
\frac{(aq/bc;q)_\infty(q^2,aq,q/a,aq^2/b^2,aq^2/c^2;q^2)_\infty}
{(aq/b,aq/c,q/b,q/c,-aq/bc;q)_\infty},
\end{equation}
valid for $|aq/bc|<1$.

Next, we give another variant of a bilateral generating function for the
bilateral $q$-ultraspherical functions. 
\begin{theorem}\label{thmgf2}
We have
\begin{align}\label{gfbil2}
\sum_{n=-\infty}^\infty (1-\ba\ga q^{n}) C_n(x;\ba,\ga\,|\,q)t^n&=
(1-\ga/\ba)(1-\ga/\ba t^2)\notag\\*
&\quad\;\times
\frac{(q,q\ga/\ba;q)_\infty^2}{(q\ga,q/\ba;q)_\infty^2}
\frac{(\ba te^{\pm\ri\ta},qe^{\pm\ri\ta}/\ba t;q)_\infty}
{(te^{\pm\ri\ta},\ga e^{\pm\ri\ta}/\ba t;q)_\infty},
\end{align}
where $|\ga/\ba|<|te^{\pm\ri\ta}|<1$.
\end{theorem}
\begin{proof}
We compute
\begin{align*}
&\sum_{n=-\infty}^\infty (1-\ba\ga q^{n}) C_n(x;\ba,\ga\,|\,q)t^n\\
&=\sum_{n=-\infty}^\infty C_n(x;\ba,\ga\,|\,q)t^n
-\ba\ga
\sum_{n=-\infty}^\infty C_n(x;\ba,\ga\,|\,q)q^nt^n
\end{align*}
using \eqref{gfbil} and \eqref{gfbil-qt}.  The generating function thus evaluates to
\begin{align*}
&\left(1-\ba\ga\frac{(1-te^{\pm\ri\ta})(1-e^{\pm\ri\ta}/\ba t)}
{(1-\ba te^{\pm\ri\ta})(1-\ga e^{\pm\ri\ta}/\ba t)}
\right)\frac{(q,q\ga/\ba;q)_\infty^2}{(q\ga,q/\ba;q)_\infty^2}
\frac{(\ba te^{\pm\ri\ta},qe^{\pm\ri\ta}/\ba t;q)_\infty}
{(te^{\pm\ri\ta},q\ga e^{\pm\ri\ta}/\ba t;q)_\infty}\\
&=
\frac{(1-\ga/\ba)(1-\ga/\ba t^2)}
{(1-\ga e^{\pm\ri\ta}/\ba t)}
\frac{(q,q\ga/\ba;q)_\infty^2}{(q\ga,q/\ba;q)_\infty^2}
\frac{(\ba te^{\pm\ri\ta},qe^{\pm\ri\ta}/\ba t;q)_\infty}
{(te^{\pm\ri\ta},q\ga e^{\pm\ri\ta}/\ba t;q)_\infty},
\end{align*}
which furnishes the claim.
\end{proof}

The $\ga\to 1$ limit of Theorem~\ref{thmgf2} readily reduces to
a second generating function for the (unilateral) continuous
$q$-ultraspherical polynomials.
\begin{corollary}\label{corgf2}
We have
\begin{equation}\label{cgfbil2}
\sum_{n=0}^\infty (1-\ba q^{n}) C_n(x;\ba\,|\,q)t^n=
(1-\ba)(1-\ba t^2)
\frac{(q\ba te^{\pm\ri\ta};q)_\infty}
{(te^{\pm\ri\ta};q)_\infty},
\end{equation}
where $|te^{\pm\ri\ta}|<1$.
\end{corollary}
This generating function has appeared before,
in equivalent normalizations, in Gasper's work on $q$-orthogonal functions
\cite{Gas81}.  A further early occurrence is in the paper of Al-Salam,
Allaway, and Askey~\cite[p.~48, just above (4.10)]{SAA}; see also the recent
treatment by Chen and Liu~\cite[Eq.~(4.1)]{CL}.
It can also be recovered as a specialization of the bilinear generating
functions for continuous $q$-ultraspherical polynomials studied by Gasper and
Rahman~\cite{GRPoisson}.

\subsection{Askey--Wilson operator and Rodrigues formulae}

The action of $\mathcal D_q$ on the bilateral $q$-ultraspherical functions
is
\begin{equation}\label{DqCn}
\mathcal D_q C_n(x;\ba,\ga\,|\,q)=
\frac{2(1-\ba)^2}{\ga(1-q)(1-\ba/\ga)}q^{\frac{1-n}2}
C_{n-1}(x;q\ba,\ga\,|\,q).
\end{equation}
This follows readily from the explicit generating function
in Theorem~\ref{thmgf}.  We have
\begin{align*}
&\mathcal D_q 
\sum_{n=-\infty}^\infty C_n(x;\ba,\ga\,|\,q)t^n=
\frac{(q,q\ga/\ba;q)_\infty^2}{(q\ga,q/\ba;q)_\infty^2}
\mathcal D_q 
\frac{(\ba te^{\pm\ri\ta},qe^{\pm\ri\ta}/\ba t;q)_\infty}
{(te^{\pm\ri\ta},q\ga e^{\pm\ri\ta}/\ba t;q)_\infty}\\
&=\frac{2\,(q,q\ga/\ba;q)_\infty^2}
{(q^{\frac 12}-q^{-\frac 12})(e^{\ri\ta}-e^{-\ri\ta})
(q\ga,q/\ba;q)_\infty^2}\\*
&\quad\;\times
\bigg(
\frac{(q^{\frac 12}\ba te^{\ri\ta},q^{-\frac 12}\ba te^{-\ri\ta},
q^{\frac 32}e^{\ri\ta}/\ba t,q^{\frac 12}e^{-\ri\ta}/\ba t;q)_\infty}
{(q^{\frac 12}te^{\ri\ta},q^{-\frac 12}te^{-\ri\ta},
q^{\frac 32}\ga e^{\ri\ta}/\ba t,q^{\frac 12}\ga e^{-\ri\ta}/\ba t;q)_\infty}\\*
&\qquad\quad-
\frac{(q^{-\frac 12}\ba te^{\ri\ta},q^{\frac 12}\ba te^{-\ri\ta},
q^{\frac 12}e^{\ri\ta}/\ba t,q^{\frac 32}e^{-\ri\ta}/\ba t;q)_\infty}
{(q^{-\frac 12}te^{\ri\ta},q^{\frac 12}te^{-\ri\ta},
q^{\frac 12}\ga e^{\ri\ta}/\ba t,q^{\frac 32}\ga e^{-\ri\ta}/\ba t;q)_\infty}
\bigg)\\
&=\frac{2\,(q,q\ga/\ba;q)_\infty^2}
{(q^{\frac 12}-q^{-\frac 12})(e^{\ri\ta}-e^{-\ri\ta})
(q\ga,q/\ba;q)_\infty^2}
\frac{(q^{\frac 12}\ba te^{\pm\ri\ta},q^{\frac 32}e^{\pm\ri\ta}/\ba t;q)_\infty}
{(q^{-\frac 12}te^{\pm\ri\ta},q^{\frac 12}\ga e^{\pm\ri\ta}/\ba t;q)_\infty}\\*
&\quad\;\times
\big((1- q^{-\frac 12}\ba te^{-\ri\ta})(1-q^{\frac 12}e^{-\ri\ta}/\ba t)
(1-q^{-\frac 12}te^{\ri\ta})(1-q^{\frac 12}\ga e^{\ri\ta}/\ba t)\\*
&\qquad\quad -(1-q^{-\frac 12}\ba t e^{\ri\ta})(1-q^{\frac 12}e^{\ri\ta}/\ba t)
(1-q^{-\frac 12}te^{-\ri\ta})(1-q^{\frac 12}\ga e^{-\ri\ta}/\ba t)\big)\\
&=\frac{2\,(q,q\ga/\ba;q)_\infty^2}
{(q^{\frac 12}-q^{-\frac 12})(e^{\ri\ta}-e^{-\ri\ta})
(q\ga,q/\ba;q)_\infty^2}
\frac{(q^{\frac 12}\ba te^{\pm\ri\ta},q^{\frac 32}e^{\pm\ri\ta}/\ba t;q)_\infty}
{(q^{-\frac 12}te^{\pm\ri\ta},q^{\frac 12}\ga e^{\pm\ri\ta}/\ba t;q)_\infty}\\*
&\quad\;\times q^{-\frac 12}\ga te^{-\ri\ta} (1-\ba/\ga)(1-e^{2\ri\ta})
(1-q^{\frac 12}e^{\ri\ta}/\ba t)(1-q^{\frac 12}e^{-\ri\ta}/\ba t)\\
&=\frac{2\ga t(1-\ba/\ga)(q,q\ga/\ba;q)_\infty^2}
{(1-q)(q\ga,q/\ba;q)_\infty^2}
\frac{(q^{\frac 12}\ba te^{\pm\ri\ta},
q^{\frac 12}e^{\pm\ri\ta}/\ba t;q)_\infty}
{(q^{-\frac 12}te^{\pm\ri\ta},q^{\frac 12}\ga e^{\pm\ri\ta}/\ba t;q)_\infty}\\
&=\frac{2t(1-\ba)^2}
{\ga(1-q)(1-\ba/\ga)}
\sum_{n=-\infty}^\infty C_n(x;q\ba,\ga\,|\,q)\big(tq^{-\frac 12}\big)^n\\
&=\frac{2(1-\ba)^2}
{\ga(1-q)(1-\ba/\ga)}
\sum_{n=-\infty}^\infty  q^{\frac{1-n}2}C_{n-1}(x;q\ba,\ga\,|\,q)t^n,
\end{align*}
from which the claimed result follows from taking coefficients of $t^n$.

\begin{remark}
The formula also gives the corresponding action of any right-inverse of
$\mathcal D_q$.  If $\mathcal I_q$ denotes such a right-inverse with the
additive $q$-constant fixed to be zero, then, away from the singular parameter
values,
\begin{equation}\label{IqCn}
\mathcal I_q C_n(x;q\ba,\ga\,|\,q)=
\frac{\ga(1-q)(1-\ba/\ga)}{2(1-\ba)^2}q^{\frac n2}
C_{n+1}(x;\ba,\ga\,|\,q).
\end{equation}
Brown and Ismail~\cite{BrownIsmailRI} constructed a distinguished such
right-inverse as an integral transform on
$L^2((-1,1),(1-x^2)^{-1/2}\rd x)$; its kernel is expressed in terms of
$\vartheta_4'/\vartheta_4$.  Thus, for that particular choice of
$\mathcal I_q$, \eqref{IqCn} is an integral representation of the raising
step from $C_n(x;q\ba,\ga\,|\,q)$ to $C_{n+1}(x;\ba,\ga\,|\,q)$.
For a different choice of right-inverse, the right-hand side of
\eqref{IqCn} is changed by an element of the kernel of $\mathcal D_q$.
Equivalently, if $H(x)=\widetilde H(z)$ with $x=(z+z^{-1})/2$, the ambiguity
is a $q$-constant, $\widetilde H(qz)=\widetilde H(z)$.  In function classes
where the only such elements are constants, this reduces to an
$x$-independent constant.

Repeated use of \eqref{DqCn} gives the corresponding lowering identity, for
$n\ge0$,
\begin{equation}\label{DqRodriguesType}
\mathcal D_q^n C_n(x;\ba,\ga\,|\,q)=
\left(\frac{2}{\ga(1-q)}\right)^n
q^{-\binom n2/2}
\frac{(\ba;q)_n^2}{(\ba/\ga;q)_n}
C_0(x;q^n\ba,\ga\,|\,q).
\end{equation}
Inverting these lowering steps with the normalized right-inverse gives a
Rodrigues-type raising construction,
\begin{equation}\label{IqRodriguesType}
C_n(x;\ba,\ga\,|\,q)=
\left(\frac{2}{\ga(1-q)}\right)^n
q^{-\binom n2/2}
\frac{(\ba;q)_n^2}{(\ba/\ga;q)_n}
\mathcal I_q^n C_0(x;q^n\ba,\ga\,|\,q),
\end{equation}
again up to the same possible $q$-constant ambiguities at each integration
step.  The parameter shift in the seed is forced by \eqref{DqCn}: applying
$\mathcal D_q$ lowers the order and simultaneously sends $\ba$ to $q\ba$.
Thus \eqref{DqRodriguesType}--\eqref{IqRodriguesType} are not yet a classical
Rodrigues formula in which repeated applications of $\mathcal D_q$ act on a
modified weight, but they do give the precise operator-theoretic form of a
Rodrigues-type construction starting from the seed $C_0$.
\end{remark}

The following weight-based Rodrigues formulae use the standard lift
$x=(z+z^{-1})/2$.  Since the weights contain the factor $(1-x^2)^{-1/2}$, we
use the same symbol $\mathcal D_q$ for the lifted divided difference
\begin{equation*}
(\mathcal D_q F)(z)=
\frac{F(q^{\frac12}z)-F(q^{-\frac12}z)}
{(q^{\frac12}-q^{-\frac12})(z-z^{-1})/2}.
\end{equation*}
This agrees with the preceding definition when $F$ comes from a symmetric
function of $z$.  Define the three lifted densities
\begin{subequations}\label{bilRodWeights}
\begin{align}
\mathcal W_{\ba,\ga}^{(1)}(z)
&=\frac{2\ri}{z-z^{-1}}
\frac{(z^{\pm2},q\ga z^{\pm2}/\ba;q)_\infty}
{(\ba z^{\pm2},qz^{\pm2}/\ba;q)_\infty},
\label{bilRodWeight}\\*
\mathcal W_{\ba,\ga}^{(2)}(z)
&=\mathcal W_{1/\ga,1/\ba}^{(1)}(z)
=\frac{2\ri}{z-z^{-1}}
\frac{(z^{\pm2},q\ga z^{\pm2}/\ba;q)_\infty}
{(z^{\pm2}/\ga,q\ga z^{\pm2};q)_\infty},
\label{bilRodWeightDual}\\*
\mathcal W_{\ba,\ga}^{(4)}(z)
&=\frac{2\ri}{z-z^{-1}}
\frac{(z^{\pm2},q\ga z^{\pm2}/\ba,
z^{\pm2}/\ba\ga,q\ba\ga z^{\pm2};q)_\infty}
{(\ba z^{\pm2},qz^{\pm2}/\ba,
z^{\pm2}/\ga,q\ga z^{\pm2};q)_\infty}.
\label{bilRodWeightFour}
\end{align}
\end{subequations}
For $z=e^{\ri\ta}$, $0<\ta<\pi$, these are precisely the continuous densities
appearing in \eqref{orth-id}, \eqref{orth-id-neg}, and \eqref{orth-id-four},
respectively.  At $\ga=1$ the first lifted density reduces to the ordinary Rogers density.
The second density is its dual under
$(\ba,\ga)\mapsto(1/\ga,1/\ba)$; for this reason its Rodrigues formula is
naturally oriented toward the negative-index functions $C_{-n}$.

The following weighted backward shifts are the first-order identities behind
the three Rodrigues formulae.  They are $q$-difference identities for the
displayed densities and do not use the full orthogonality relations proved
later in Section~\ref{sec:orth}.
\begin{proposition}[Bilateral weighted backward shifts]\label{prop:bilbackshift}
Let $m\in\Z$, and suppose that the parameters avoid the poles in the displayed
expressions.  Then, with $x=(z+z^{-1})/2$,
\begin{subequations}\label{bilBackShifts}
\begin{align}
\mathcal D_q\!
\left[\mathcal W_{\ba,\ga}^{(1)}(z)C_m(x;\ba,\ga\,|\,q)\right]
&=A_m^{(1)}(\ba,\ga;q)\,
\mathcal W_{\ba/q,\ga}^{(1)}(z)C_{m+1}(x;\ba/q,\ga\,|\,q),
\label{bilBackShift}\\*
\mathcal D_q\!
\left[\mathcal W_{\ba,\ga}^{(2)}(z)C_m(x;\ba,\ga\,|\,q)\right]
&=A_m^{(2)}(\ba,\ga;q)\,
\mathcal W_{\ba,q\ga}^{(2)}(z)C_{m-1}(x;\ba,q\ga\,|\,q),
\label{bilBackShiftDual}\\*
\mathcal D_q\!
\left[\mathcal W_{\ba,\ga}^{(4)}(z)C_m(x;\ba,\ga\,|\,q)\right]
&=A_m^{(4)}(\ba,\ga;q)\,
\mathcal W_{\ba/q,\ga}^{(4)}(z)C_{m+1}(x;\ba/q,\ga\,|\,q),
\label{bilBackShiftFour}
\end{align}
\end{subequations}
where
\begin{subequations}\label{bilBackShiftCoeffs}
\begin{align}
A_m^{(1)}(\ba,\ga;q)
&=\frac{2q^{1-m/2}(\ba-q\ga)(1-\ga^2q^{m+1})(1-\ba^2q^{m-1})}
{(1-q)(\ba-q)^2},
\label{bilBackShiftCoeff}\\*
A_m^{(2)}(\ba,\ga;q)
&=\frac{2q^{3m/2-1}(\ba-q\ga)(1-q^{1-m}/\ba^2)(1-q^{-m-1}/\ga^2)}
{(1-q)(1-q\ga)^2},
\label{bilBackShiftCoeffDual}\\*
A_m^{(4)}(\ba,\ga;q)
&=\frac{A_m^{(1)}(\ba,\ga;q)}{\ba^2\ga^2}.
\label{bilBackShiftCoeffFour}
\end{align}
\end{subequations}
\end{proposition}
\begin{proof}
The identities are first verified in a domain of absolute convergence and then
continued meromorphically in the parameters.  For the first one, the product
part of $\mathcal W_{\ba,\ga}^{(1)}(q^{1/2}z)$ can be rewritten relative to
$\mathcal W_{\ba/q,\ga}^{(1)}(z)$ by using
$(qa;q)_\infty=(a;q)_\infty/(1-a)$ and
$(q^{-1}a;q)_\infty=(1-q^{-1}a)(a;q)_\infty$.  This gives
\begin{align}\label{bilRodRatioPlus}
\frac{\mathcal W_{\ba,\ga}^{(1)}(q^{1/2}z)}
{\mathcal W_{\ba/q,\ga}^{(1)}(z)}
&=\frac{(z-z^{-1})}{(q^{1/2}z-q^{-1/2}z^{-1})}
\frac{(1-q^{-1}z^{-2})}{(1-z^2)}\notag\\*
&\quad\times
\frac{(1-\ga z^{-2}/\ba)(1-q\ga z^{-2}/\ba)
(1-\ba z^2/q)(1-\ba z^2)}
{(1-z^{-2}/\ba)(1-qz^{-2}/\ba)}.
\end{align}
The analogous formula for $\mathcal W_{\ba,\ga}^{(1)}(q^{-1/2}z)$ is obtained
from \eqref{bilRodRatioPlus} by the involution $z\mapsto z^{-1}$.
Insert the bilateral series \eqref{defbilC} into the left-hand side of
\eqref{bilBackShift}, use the two ratio formulae, and rewrite the result with
the common factor $\mathcal W_{\ba/q,\ga}^{(1)}(z)$.  The coefficient of
$z^{m+1-2k}$ is
\[
A_m^{(1)}(\ba,\ga;q)
\frac{(\ba/q;q)_k(\ba/q;q)_{m+1-k}}
{(q\ga;q)_k(q\ga;q)_{m+1-k}},
\]
which is the corresponding coefficient in the right-hand side of
\eqref{bilBackShift}.

For the dual identity \eqref{bilBackShiftDual}, use
$\mathcal W_{\ba,\ga}^{(2)}=\mathcal W_{1/\ga,1/\ba}^{(1)}$ together with the
symmetry \eqref{sym}.  Applying \eqref{bilBackShift} with parameters
$(1/\ga,1/\ba)$ and order $-m$, and then converting back by \eqref{sym}, gives
\eqref{bilBackShiftDual}.  The scalar prefactor is exactly the identity
\[
\left(\frac{q\ga}{\ba}\right)^{-m}
\left(\frac{q^2\ga}{\ba}\right)^{m-1}
 A_{-m}^{(1)}(1/\ga,1/\ba;q)=A_m^{(2)}(\ba,\ga;q).
\]
For \eqref{bilBackShiftFour}, write
\[
\mathcal W_{\ba,\ga}^{(4)}(z)=\mathcal W_{\ba,\ga}^{(1)}(z)
\frac{(z^{\pm2}/\ba\ga,q\ba\ga z^{\pm2};q)_\infty}
{(z^{\pm2}/\ga,q\ga z^{\pm2};q)_\infty}.
\]
The same coefficient comparison as above, now relative to the common factor
$\mathcal W_{\ba/q,\ga}^{(4)}(z)$, shows that the extra product multiplies
the coefficient in \eqref{bilBackShiftCoeff} by $(\ba\ga)^{-2}$.  This proves
\eqref{bilBackShiftFour}.
\end{proof}

\begin{remark}
For $\ga=1$, the coefficient in \eqref{bilBackShiftCoeff} becomes
$A_m^{(1)}(\ba,1;q)=
-2q^{-m/2}(1-q^{m+1})(1-\ba^2q^{m-1})/
((1-q)(1-\ba q^{-1}))$, so \eqref{bilBackShift} is the standard Rogers weighted backward shift
in the present normalization.
\end{remark}

Iterating Proposition~\ref{prop:bilbackshift} gives the three
Rodrigues-type formulae.
\begin{theorem}[Bilateral Rodrigues-type formulae]\label{thm:bilrodrigues}
Let $n=0,1,2,\ldots$, and suppose that the parameters avoid the poles in the
displayed expressions.  Then, with $x=(z+z^{-1})/2$,
\begin{subequations}\label{bilRodriguesAll}
\begin{align}
\mathcal W_{\ba,\ga}^{(1)}(z)C_n(x;\ba,\ga\,|\,q)
&=B_n^{(1)}(\ba,\ga;q)\,
\mathcal D_q^n\!
\left[\mathcal W_{\ba q^n,\ga}^{(1)}(z)C_0(x;\ba q^n,\ga\,|\,q)\right],
\label{bilRodrigues}\\
\mathcal W_{\ba,\ga}^{(2)}(z)C_{-n}(x;\ba,\ga\,|\,q)
&=B_n^{(2)}(\ba,\ga;q)\,
\mathcal D_q^n\!
\left[\mathcal W_{\ba,\ga q^{-n}}^{(2)}(z)C_0(x;\ba,\ga q^{-n}\,|\,q)\right],
\label{bilRodriguesDual}\\
\mathcal W_{\ba,\ga}^{(4)}(z)C_n(x;\ba,\ga\,|\,q)
&=B_n^{(4)}(\ba,\ga;q)\,
\mathcal D_q^n\!
\left[\mathcal W_{\ba q^n,\ga}^{(4)}(z)C_0(x;\ba q^n,\ga\,|\,q)\right],
\label{bilRodriguesFour}
\end{align}
\end{subequations}
where
\begin{subequations}\label{bilRodriguesCoeffs}
\begin{align}
B_n^{(1)}(\ba,\ga;q)
&=\left(\frac{q-1}{2}\right)^n q^{n(n-1)/4}
\frac{(\ba;q)_n^2}
{\ga^n(\ba/\ga;q)_n(\ga^2q;q)_n(\ba^2q^n;q)_n},
\label{bilRodriguesCoeff}\\*
B_n^{(2)}(\ba,\ga;q)
&=\left(\frac{1-q}{2}\right)^n q^{(3n^2+n)/4}
\frac{(\ga q^{1-n};q)_n^2}
{\ba^n(\ga q^{1-n}/\ba;q)_n(q/\ba^2;q)_n(q^n/\ga^2;q)_n},
\label{bilRodriguesCoeffDual}\\*
B_n^{(4)}(\ba,\ga;q)
&=\ba^{2n}\ga^{2n}q^{n(n+1)}B_n^{(1)}(\ba,\ga;q)
\label{bilRodriguesCoeffFourFactor}\\*
&=\left(\frac{q-1}{2}\right)^n q^{(5n^2+3n)/4}
\frac{\ba^{2n}\ga^n(\ba;q)_n^2}
{(\ba/\ga;q)_n(\ga^2q;q)_n(\ba^2q^n;q)_n}.
\label{bilRodriguesCoeffFour}
\end{align}
\end{subequations}
\end{theorem}
\begin{proof}
Put $\ba_j=\ba q^j$.  Starting from
$\mathcal W_{\ba_n,\ga}^{(1)}(z)C_0(x;\ba_n,\ga\,|\,q)$, apply
\eqref{bilBackShift} successively with
$(m,\ba)=(0,\ba_n),(1,\ba_{n-1}),\ldots,(n-1,\ba_1)$.  This gives the same identity as
\eqref{bilRodrigues}, with coefficient
$P_n^{(1)}(\ba,\ga;q)=\prod_{r=0}^{n-1}A_r^{(1)}(\ba q^{n-r},\ga;q)$ on the right-hand side before inversion.  Using \eqref{bilBackShiftCoeff}, we have
\begin{align*}
A_r^{(1)}(\ba q^{n-r},\ga;q)
&=\frac{2q^{1-r/2}}{(1-q)}
\frac{(\ba q^{n-r}-q\ga)(1-\ga^2q^{r+1})(1-\ba^2q^{2n-r-1})}
{(\ba q^{n-r}-q)^2}\\*
&=\frac{2q^{-r/2}}{(1-q)}
\frac{(\ba q^{n-r-1}-\ga)(1-\ga^2q^{r+1})(1-\ba^2q^{2n-r-1})}
{(1-\ba q^{n-r-1})^2}.
\end{align*}
Multiplication over $r=0,\ldots,n-1$ gives
\[
P_n^{(1)}(\ba,\ga;q)=
\left(\frac{2}{1-q}\right)^n q^{-n(n-1)/4}
\frac{\prod_{s=0}^{n-1}(\ba q^s-\ga)(\ga^2q;q)_n(\ba^2q^n;q)_n}
{(\ba;q)_n^2}.
\]
Thus $B_n^{(1)}=(P_n^{(1)})^{-1}$.  Rewriting
$\ga^n(\ba/\ga;q)_n=\prod_{s=0}^{n-1}(\ga-\ba q^s)=(-1)^n\prod_{s=0}^{n-1}(\ba q^s-\ga)$ gives \eqref{bilRodriguesCoeff}.

For the dual formula, start from
$\mathcal W_{\ba,\ga q^{-n}}^{(2)}(z)C_0(x;\ba,\ga q^{-n}\,|\,q)$ and apply
\eqref{bilBackShiftDual} successively with
$(m,\ga)=(0,\ga q^{-n}),(-1,\ga q^{-n+1}),\ldots,(-(n-1),\ga q^{-1})$.  The product of the coefficients is
\begin{align*}
P_n^{(2)}(\ba,\ga;q)
&=\prod_{r=0}^{n-1}A_{-r}^{(2)}(\ba,\ga q^{-n+r};q)\\*
&=\left(\frac{2}{1-q}\right)^nq^{-(3n^2+n)/4}
\frac{\ba^n(\ga q^{1-n}/\ba;q)_n(q/\ba^2;q)_n(q^n/\ga^2;q)_n}
{(\ga q^{1-n};q)_n^2}.
\end{align*}
Hence $B_n^{(2)}=(P_n^{(2)})^{-1}$, giving
\eqref{bilRodriguesCoeffDual} and \eqref{bilRodriguesDual}.

Finally, the four-factor formula is obtained by the same iteration as the
first formula, but with \eqref{bilBackShiftFour}.  By
\eqref{bilBackShiftCoeffFour},
$P_n^{(4)}=P_n^{(1)}/(\ba^{2n}\ga^{2n}q^{n(n+1)})$, hence
$B_n^{(4)}=\ba^{2n}\ga^{2n}q^{n(n+1)}B_n^{(1)}$.  This proves \eqref{bilRodriguesCoeffFourFactor}; substituting
\eqref{bilRodriguesCoeff} gives \eqref{bilRodriguesCoeffFour}.
\end{proof}

When $\ga=1$, the seed in the first formula \eqref{bilRodrigues} collapses to
$C_0(x;\ba q^n,1\,|\,q)=1$ and $\mathcal W_{\ba,1}^{(1)}$ is the ordinary
Rogers density.  Hence \eqref{bilRodrigues} reduces to the standard
Rodrigues-type formula for the continuous $q$-ultraspherical polynomials,
\begin{equation}\label{classicalRodriguesFromBilateral}
\mathcal W_{\ba,1}^{(1)}(z)C_n(x;\ba\,|\,q)
=\left(\frac{q-1}{2}\right)^n q^{n(n-1)/4}
\frac{(\ba;q)_n}{(q;q)_n(\ba^2q^n;q)_n}
\mathcal D_q^n\mathcal W_{\ba q^n,1}^{(1)}(z).
\end{equation}
This is the Rogers formula in the present normalization; compare
\cite[Eq.~(3.10.24)]{KS}.

\begin{remark}
In the genuinely bilateral case, the seeds in \eqref{bilRodrigues},
\eqref{bilRodriguesDual}, and \eqref{bilRodriguesFour} cannot be replaced by
the shifted densities alone.  Indeed $C_0(x;\ba,\ga\,|\,q)$ is the
nontrivial bilateral series displayed in \eqref{C0spec}; it is exactly the
information left after the corresponding lowering steps.  For terminating
parameters, such as $\ga=q^s$ with $s\ge0$, this seed becomes a finite
Laurent-polynomial factor, whereas for generic $\ga$ it remains a
nonterminating bilateral contribution.
\end{remark}

\subsection{A special integral}

In the ordinary polynomial case $C_0$ is a constant, whereas in the bilateral
case $C_0(x;\ba,\ga\,|\,q)$ is generally a nontrivial function, as in
\eqref{C0spec}.  The following special integral shows that this seed still has
a simple orthogonality-type evaluation when $\ga=1/\ba$:
\begin{equation}\label{almostorth}
\frac 1{2\pi}\int_{-1}^1C_n(x;\ba,1/\ba\,|\,q)\,
\frac{(e^{\pm 2\ri\ta};q)_\infty}{(\ba e^{\pm 2\ri\ta};q)_\infty}
\frac{\rd x}{\sqrt{1-x^2}}=
\frac{(q;q)_\infty^2(\ba,q/\ba^2;q)_\infty}{(q/\ba;q)_\infty^3(\ba^2;q)_\infty}\,
\da_{n,0}.
\end{equation}
One can prove \eqref{almostorth} by employing the integral evaluation
in \eqref{thm1322}.  Indeed, we have
\begin{align}\label{comp}
\frac 1{2\pi}\int_{-1}^1\sum_{n=-\infty}^\infty
C_n(x;\ba,1/\ba\,|\,q)t^n\,
\frac{(e^{\pm 2\ri\ta};q)_\infty}{(\ba e^{\pm 2\ri\ta};q)_\infty}
\frac{\rd x}{\sqrt{1-x^2}}&\notag\\=
\frac 1{2\pi}\int_{-1}^1
\frac{(q,q/\ba^2;q)_\infty^2}{(q/\ba;q)_\infty^4}
\frac{(\ba te^{\pm\ri\ta},qe^{\pm\ri\ta}/\ba t,e^{\pm 2\ri\ta};q)_\infty}
{(te^{\pm\ri\ta},qe^{\pm\ri\ta}/\ba^2 t,\ba e^{\pm 2\ri\ta};q)_\infty}
\frac{\rd x}{\sqrt{1-x^2}}&\notag\\
=\frac{(q,q/\ba^2;q)_\infty^2}{(q/\ba;q)_\infty^4}\,
\frac{(\ba,q\ba;q)_\infty}{(q,\ba^2;q)_\infty}\,
{}_2\phi_1\!\left[\begin{matrix}
\ba^2,\ba\\q\ba\end{matrix};q,q/\ba^2\right]&\notag\\
=\frac{(q,q/\ba^2;q)_\infty^2}{(q/\ba;q)_\infty^4}\,
\frac{(\ba,q\ba;q)_\infty}{(q,\ba^2;q)_\infty}\,
\frac{(q/\ba,q;q)_\infty}{(q\ba,q/\ba^2;q)_\infty}&\notag\\=
\frac{(q;q)_\infty^2(\ba,q/\ba^2;q)_\infty}{(q/\ba;q)_\infty^3(\ba^2;q)_\infty}&.
\end{align}
In the first equality we use \eqref{gfbil}, in the
second equality \eqref{thm1322}, and
in the third we use the $q$-Gau{\ss} summation (cf.\
\cite[Equation (II.8)]{GR}),
\begin{equation}
{}_2\phi_1\!\left[\begin{matrix}
a,b\\c\end{matrix};q,\frac c{ab}\right]=
\frac{(c/a,c/b;q)_\infty}{(c,c/ab;q)_\infty},\qquad |c/ab|<1.
\end{equation}
Comparison of coefficients of $t^n$ on the left- and
right-hand sides of \eqref{comp}
establishes \eqref{almostorth}.

\section{Asymptotics}\label{sec:asym}

\subsection{Asymptotics on the interval}

Throughout this section $0<|q|<1$.  We first consider real $x$ with
$-1<x<1$.  In the cosine forms below the parameters $q$, $\ba$, and $\ga$ are
taken real, away from poles, so that the two exponential terms are complex
conjugates; for complex parameters the corresponding two-exponential forms
should be used instead.  The leading term, obtained from \eqref{2terma}, is
\begin{subequations}\label{asymld}
\begin{align}
C_n(x;\ba,\ga\,|\,q)&\asymp\frac{(q,\ba,q\ga/\ba;q)_\infty}
{(q\ga,q\ga,q/\ba;q)_\infty}\notag\\*
&\quad\;\times\left(\frac{(\ba e^{-2\ri\ta},q e^{2\ri\ta}/\ba;q)_\infty}
{(e^{-2\ri\ta},q\ga e^{2\ri\ta}/\ba;q)_\infty}e^{\ri n\ta}+
\frac{(\ba e^{2\ri\ta},q e^{-2\ri\ta}/\ba;q)_\infty}
{(e^{2\ri\ta},q\ga e^{-2\ri\ta}/\ba;q)_\infty}e^{-\ri n\ta}
\right)\notag\\
&=2\frac{(q,\ba,q\ga/\ba;q)_\infty}
{(q\ga,q\ga,q/\ba;q)_\infty}|A_{\ba,\ga}(e^{\ri\ta})|
\cos(n\ta-\alpha),
\end{align}
where
\begin{equation}
A_{\ba,\ga}(z)=
\frac{(\ba z^2,q/\ba z^2;q)_\infty}
{(z^2,q\ga/\ba z^2;q)_\infty}\quad\text{and}\quad \alpha=\arg\!\big(A_{\ba,\ga}(e^{\ri\ta})\big).
\end{equation}
\end{subequations}
For $\ga\to 1$ \eqref{asymld} reduces to \cite[Eq.~(3.11)]{AI},
derived there by Darboux's method.
The same representation gives a complete finite asymptotic expansion.  For
$0<\ta<\pi$ and any fixed positive integer $N$, one obtains
\begin{subequations}\label{asymfe}
\begin{align}
&C_n(x;\ba,\ga\,|\,q)\notag\\*
&\asymp
2\frac{(q,\ba,q\ga/\ba;q)_\infty}
{(q\ga,q\ga,q/\ba;q)_\infty}
|A_{\ba,\ga}(e^{\ri\ta})|
\cos(n\ta-\alpha)\notag\\*
&\quad\;+2\frac{(q,\ba,q\ga/\ba;q)_\infty}
{(q\ga,q\ga,q/\ba;q)_\infty}\notag\\*
&\qquad\;\;\times
\sum_{k=1}^N\bigg(\frac{(q\ga/\ba;q)_k}{(q;q)_k}
|D_{k;\ba,\ga}(e^{\ri\ta})|
\cos\!\big(n\ta-\phi_k\big) \ba^{2k}q^{nk}\bigg)
+O(q^{(N+1)n}),
\end{align}
where
\begin{equation}
D_{k;\ba,\ga}(z)=
\frac{(\ba z^2,q/\ba z^2;q)_\infty}
{(z^2,q\ga/\ba z^2;q)_\infty}\frac{(q\ga/\ba z^2;q)_k}
{(q/z^2;q)_k}
\quad\text{for $k=1,2,\ldots$},
\end{equation}
and
\begin{equation}
\phi_k=\arg\!\big(D_{k;\ba,\ga}(e^{\ri\ta})\big).
\end{equation}
\end{subequations}
The transformed representation \eqref{2termb} gives the equivalent expansion
\begin{subequations}\label{asymfe2}
\begin{align}
&C_n(x;\ba,\ga\,|\,q)\notag\\*
&\asymp
2\frac{(q,\ba,q\ga/\ba;q)_\infty}
{(q\ga,q\ga,q/\ba;q)_\infty}
|A_{\ba,\ga}(e^{\ri\ta})|
\cos(n\ta-\alpha)\notag\\*
&\quad\;+2\frac{(q,\ba,q\ga/\ba;q)_\infty}
{(q\ga,q\ga,q/\ba;q)_\infty}\notag\\
&\qquad\;\;\times
\sum_{k=1}^N\bigg(\frac{(\ba/\ga;q)_k}{(q;q)_k}
|E_{k;\ba,\ga}(e^{\ri\ta})|
\cos\!\big(n\ta-\psi_k\big) \ga^{2k}q^{nk+k}\bigg)
+O(q^{(N+1)(n+1)}),
\end{align}
where
\begin{equation}
E_{k;\ba,\ga}(z)=
\frac{(\ba z^2,q/\ba z^2;q)_\infty}
{(z^2,q\ga/\ba z^2;q)_\infty}\frac{(\ba/\ga z^2;q)_k}
{(q/z^2;q)_k}
\quad\text{for $k=1,2,\ldots$},
\end{equation}
and
\begin{equation}
\psi_k=\arg\!\big(E_{k;\ba,\ga}(e^{\ri\ta})\big).
\end{equation}
\end{subequations}
For $\ga\to1$, \eqref{asymfe2} reduces to \cite[Eq.~(5.3)]{RV}.

\subsection{Off-interval and negative-order asymptotics}

If $x\in\C\setminus[-1,1]$, \eqref{2term} can be written as
\begin{align*}
C_n(x;\ba,\ga\,|\,q)
&=\frac{(q,\ba,q\ga/\ba;q)_\infty}
{(q\ga,q\ga,q/\ba;q)_\infty}\notag\\*
&\quad\;\times\!\left(\frac{(\ba/z^2,q z^2/\ba;q)_\infty}
{(1/z^2,q\ga z^2/\ba;q)_\infty}z^n
{}_2\phi_1\!\left[\begin{matrix}q\ga/\ba,q\ga z^2/\ba\\
q z^2\end{matrix};q,\ba^2q^n\right]
\right.\notag\\*
&\qquad\quad +\left.
\frac{(\ba z^2,q/\ba z^2;q)_\infty}
{(z^2,q\ga/\ba z^2;q)_\infty}z^{-n}
{}_2\phi_1\!\left[\begin{matrix}q\ga/\ba,q\ga/\ba z^2\\
q/z^2\end{matrix};q,\ba^2q^n\right]
\right)\\
&=\frac{(q,\ba,q\ga/\ba,q\ga^2;q)_\infty}
{(q\ga,q\ga,q/\ba,\ba^2;q)_\infty}
\frac{(\ba^2;q)_n}{(q\ga^2;q)_n}\notag\\*
&\quad\;\times\!\left(\frac{(\ba/z^2,q z^2/\ba;q)_\infty}
{(1/z^2,q\ga z^2/\ba;q)_\infty}z^n
{}_2\phi_1\!\left[\begin{matrix}\ba/\ga,\ba z^2/\ga\\
q z^2\end{matrix};q,\ga^2q^{n+1}\right]
\right.\notag\\*
&\qquad\quad +\left.
\frac{(\ba z^2,q/\ba z^2;q)_\infty}
{(z^2,q\ga/\ba z^2;q)_\infty}z^{-n}
{}_2\phi_1\!\left[\begin{matrix}\ba/\ga,\ba/\ga z^2\\
q/z^2\end{matrix};q,\ga^2q^{n+1}\right]
\right),
\end{align*}
where $x=(z+1/z)/2$.  Since $x\in\C\setminus[-1,1]$,
we must have $|z|\neq 1$.  For given $x$, the equation $x=(z+1/z)/2$
has two solutions
$z_1$ and $z_2=1/z_1$, which we order according to $|z_1|<1<|z_2|$.
The asymptotics of $C_n(x;\ba,\ga\,|\,q)$ are then determined by the
dominant terms in the above series transformations, and are given by
either of the following two formulae:
\begin{subequations}
\begin{align}
C_n(x;\ba,\ga\,|\,q)&\asymp
\frac{(q,\ba,q\ga/\ba;q)_\infty}
{(q\ga,q\ga,q/\ba;q)_\infty}\frac{(\ba z_1^2,q z_2^2/\ba;q)_\infty}
{(z_1^2,q\ga z_2^2/\ba;q)_\infty}z_2^n\notag\\*
&\quad\;\times
\bigg(1+\sum_{k=1}^N
\frac{(q\ga/\ba,q\ga z_2^2/\ba;q)_k}{(q,qz_2^2;q)_k}
\ba^{2k}q^{nk}\bigg)+O\big(z_2^n q^{(N+1)n}\big)\\*
&\asymp\frac{(q,\ba,q\ga/\ba;q)_\infty}
{(q\ga,q\ga,q/\ba;q)_\infty}\frac{(\ba z_1^2,q z_2^2/\ba;q)_\infty}
{(z_1^2,q\ga z_2^2/\ba;q)_\infty}z_2^n\notag\\*
&\quad\;\times
\bigg(1+\sum_{k=1}^N
\frac{(\ba/\ga,\ba z_2^2/\ga;q)_k}{(q,qz_2^2;q)_k}
\ga^{2k}q^{nk+k}\bigg)
+O\big(z_2^n q^{(N+1)(n+1)}\big)\label{asymc2}
\end{align}
\end{subequations}
for $N=1,2,\ldots$.  For $\ga\to 1$ \eqref{asymc2} reduces
to \cite[Eq.~(5.7)]{RV}.

The large negative-order asymptotics follow directly from the symmetry
\eqref{sym}, by replacing $(\ba,\ga,n)$ with $(1/\ga,1/\ba,-n)$ in the
large-positive-order formulae.

\section{Shifted orthogonality of the bilateral \texorpdfstring{$q$}{q}-ultraspherical functions}
\label{sec:shifted}

Before proving the full orthogonality relations in Section~\ref{sec:orth},
we record two shifted orthogonality relations.  They are not substitutes for
full pairwise orthogonality; rather, they are useful companion identities for
the whole two-sided sequence $(C_k(x;\ba,\ga\,|\,q))_{k\in\mathbb Z}$.  Similar
``shifted'' conditions occur for wavelet matrices (cf.\ \cite{EL,KT,TK}).  In
contrast with the full orthogonality relations below, the identities in this
section do not reduce, when $\ga\to1$, to the ordinary Rogers orthogonality.
They are included because they are direct consequences of the bilateral
generating functions and because they illustrate how naturally the entire
integer-indexed family enters the theory.

Recall that the weight function $w(x\,|\,\ba)$ is defined in \eqref{orthw}; in
this section it is used with $\ba$ replaced by $\ba/\ga$.
\begin{proposition}
\label{prop:shifted}
Let $|q\ga/\ba^2|<|t|<|1/\ba|$.
Then we have
\begin{align}
\frac 1{2\pi}\int_{-1}^1
\sum_{k=-\infty}^\infty C_{m+k}(x;\ba,\ga\,|\,q)\,C_{n+k}(x;\ba,\ga\,|\,q)\,
\left(\frac{q\ga}{\ba^2}\right)^k
w(x\,|\,\ba/\ga)\rd x\notag\\
=\frac{(q;q)_\infty^3(q\ga/\ba;q)_\infty^4(\ba/\ga,q\ba/\ga;q)_\infty}{(q\ga,q/\ba;q)_\infty^4
(\ba^2/\ga^2;q)_\infty}\,{}_2\phi_1\!\left[\begin{matrix}
\ba^2/\ga^2,\ba/\ga\\q\ba/\ga\end{matrix};q,\frac {q\ga}{\ba^2}\right]
\left(\frac{\ba^2}{q\ga}\right)^n\,
\da_{m,n}&.
\end{align}
\end{proposition}
\begin{proof}
Consider the double bilateral generating function
  (using \eqref{gfbil} twice)
\begin{align*}
&\sum_{m,k=-\infty}^\infty
C_{m}(x;\ba,\ga\,|\,q)\,C_{k}(x;\ba,\ga\,|\,q)\;t_1^m\, t_2^k\\&=
\frac{(q,q\ga/\ba;q)_\infty^4}{(q\ga,q/\ba;q)_\infty^4}
\frac{(\ba t_1e^{\pm\ri\ta},qe^{\pm\ri\ta}/\ba t_1,
\ba t_2e^{\pm\ri\ta},qe^{\pm\ri\ta}/\ba t_2;q)_\infty}
{(t_1e^{\pm\ri\ta},q\ga e^{\pm\ri\ta}/\ba t_1,
t_2e^{\pm\ri\ta},q\ga e^{\pm\ri\ta}/\ba t_2;q)_\infty},
\end{align*}
and take $(t_1,t_2)=(t,q\ga/\ba^2 t)$.  This gives
\begin{align*}
&\sum_{m,k=-\infty}^\infty
C_{m}(x;\ba,\ga\,|\,q)\,C_{k}(x;\ba,\ga\,|\,q)\;t^{m-k}\,
\Big(\frac {q\ga}{\ba^2}\Big)^k\\&=
\frac{(q,q\ga/\ba;q)_\infty^4}{(q\ga,q/\ba;q)_\infty^4}
\frac{(qe^{\pm\ri\ta}/\ba t,\ba t e^{\pm\ri\ta}/\ga;q)_\infty}
{(t e^{\pm\ri\ta}, q\ga e^{\pm\ri\ta}/\ba^2 t;q)_\infty}.
\end{align*}
Integration over $x$ from $-1$ to $1$ with respect to the
$q$-ultraspherical weight function $ w(x\,|\,\ba/\ga)$ divided by $2\pi$ gives,
by the integral evaluation in \eqref{thm1322},
\begin{equation*}
\frac{(q,q\ga/\ba;q)_\infty^4}{(q\ga,q/\ba;q)_\infty^4}\frac
{(\ba/\ga,q\ba/\ga;q)_\infty}{(q,\ba^2/\ga^2;q)_\infty}
{}_2\phi_1\!\left[\begin{matrix}
\ba^2/\ga^2,\ba/\ga\\q\ba/\ga\end{matrix};q,\frac{q\ga}{\ba^2}\right]
\end{equation*}
which is independent of $t$.  Now shift indices $(m,k)$ by
$(m+k,n+k)$ and compare coefficients of $t$.
\end{proof}

A slight alteration of the proof gives the following variant of the
shifted orthogonality relation for the bilateral $q$-ultraspherical functions,
in which the product side completely factorizes.

\begin{proposition}
\label{prop:shifted2}
Let $|\ga/\ba^2|<|t|<|1/\ba|$.
Then we have
\begin{align}
\frac 1{2\pi}\int_{-1}^1
  \sum_{k=-\infty}^\infty \bigg(&C_{m+k}(x;\ba,\ga\,|\,q)\,
                                  C_{n+k}(x;\ba,\ga\,|\,q)\notag\\*
&\times  (1-\ba\ga q^{m+k})(1-\ba\ga q^{n+k})
\left(\frac{\ga}{\ba^2}\right)^k\bigg)\,
w(x\,|\,\ba/\ga)\rd x\notag\\
 & =(1-\ga/\ba)^2\frac{(q;q)_\infty^3(\ga/\ba,q\ga/\ba;q)_\infty^2
  (\ba/\ga,q\ba/\ga,q/\ga;q)_\infty}{(q\ga,q/\ba;q)_\infty^4
(\ba^2/\ga^2,\ga/\ba^2;q)_\infty}\,\ba(1+\ba)\,
\da_{m,n}.
\end{align}
\end{proposition}
\begin{proof}
  Consider the double bilateral generating function
  (using \eqref{gfbil2} twice)
\begin{align*}
&\sum_{m,k=-\infty}^\infty
  C_{m}(x;\ba,\ga\,|\,q)\,C_{k}(x;\ba,\ga\,|\,q)\,(1-\ba\ga q^m)(1-\ba\ga q^k)\;t_1^m\, t_2^k\\&=
  (1-\ga/\ba)^2(1-\ga/\ba t_1^2)(1-\ga/\ba t_2^2)\\*
  &\quad\;\times
\frac{(q,q\ga/\ba;q)_\infty^4}{(q\ga,q/\ba;q)_\infty^4}
\frac{(\ba t_1e^{\pm\ri\ta},qe^{\pm\ri\ta}/\ba t_1,
\ba t_2e^{\pm\ri\ta},qe^{\pm\ri\ta}/\ba t_2;q)_\infty}
{(t_1e^{\pm\ri\ta},\ga e^{\pm\ri\ta}/\ba t_1,
t_2e^{\pm\ri\ta},\ga e^{\pm\ri\ta}/\ba t_2;q)_\infty},
\end{align*}
and take $(t_1,t_2)=(t,\ga/\ba^2 t)$.  This gives
\begin{align*}
&\sum_{m,k=-\infty}^\infty
  C_{m}(x;\ba,\ga\,|\,q)\,C_{k}(x;\ba,\ga\,|\,q)\,
  (1-\ba\ga q^m)(1-\ba\ga q^k)\;t^{m-k}\,
\Big(\frac {\ga}{\ba^2}\Big)^k\\&=
  (1-\ga/\ba)^2(1-\ga/\ba t^2)(1-\ba^3t^2/\ga)
  \frac{(q,q\ga/\ba;q)_\infty^4}{(q\ga,q/\ba;q)_\infty^4}
\frac{(qe^{\pm\ri\ta}/\ba t,q\ba t e^{\pm\ri\ta}/\ga;q)_\infty}
{(t e^{\pm\ri\ta}, \ga e^{\pm\ri\ta}/\ba^2 t;q)_\infty}.
\end{align*}
Integration over $x$ from $-1$ to $1$ with respect to the
$q$-ultraspherical weight function $ w(x\,|\,\ba/\ga)$ divided by $2\pi$ gives,
by the integral evaluation in \eqref{orth-idgf3-cont},
\begin{equation*}
 (1-\ga/\ba)^2 \frac{(q,q\ga/\ba;q)_\infty^4}{(q\ga,q/\ba;q)_\infty^4}
  \frac{(\ga/\ba,q\ga/\ba,q/\ga;q)_\infty}
{(q,\ba^2/\ga^2,\ga/\ba^2;q)_\infty}\,\ba(1+\ba)
\end{equation*}
which is independent of $t$.  Now shift indices $(m,k)$ by
$(m+k,n+k)$ and compare coefficients of $t$.
\end{proof}

\section{Full orthogonality of the bilateral \texorpdfstring{$q$}{q}-ultraspherical functions}\label{sec:orth}

\subsection{Mass aggregates and two-factor orthogonality}

This section contains the main orthogonality results of the paper.  The point is
not merely to continue the integral part of the ordinary Rogers orthogonality,
but to identify the full continued orthogonality functionals, including the
residue mass aggregates that replace the moving finite mass points.

The analytic continuation step used below is close in spirit to Ismail's
argument for bilateral summations: prove an identity on a terminating or otherwise
specializing set with an accumulation point, and then extend it by analyticity
or meromorphic continuation.  This argument was used by Ismail for Ramanujan's
${}_1\psi_1$ summation~\cite{I} and by Askey and Ismail for Bailey's
very-well-poised ${}_6\psi_6$ summation~\cite{AI6}.  Here the same idea is
applied not to a single bilateral series, but to the integral-plus-residue
orthogonality functional.

The pure integral relations suggested by formal analytic continuation are
incomplete in general.  In the terminating specializations $\ga=q^s$ the
bilateral functions reduce to ordinary continuous $q$-ultraspherical
polynomials with shifted parameter $\ba q^{-s}$; once this shifted parameter
passes the pure integral range, the ordinary Askey--Wilson mass points have
to be included.  For generic $\ga$ these moving finite mass sums are replaced
by the following analytic residue aggregates.

For $a\ne0$, put
\begin{equation*}
 x_a=\frac{a^{1/2}+a^{-1/2}}2
\end{equation*}
and define the parity-symmetrized product
\begin{equation}\label{Cmndef}
\mathcal C_{m,n}(a)=\frac12\Big(
C_m(x_a;\ba,\ga\,|\,q)C_n(x_a;\ba,\ga\,|\,q)
+C_m(-x_a;\ba,\ga\,|\,q)C_n(-x_a;\ba,\ga\,|\,q)
\Big).
\end{equation}
The two analytic mass aggregates needed below are
\begin{subequations}\label{massaggs}
\begin{align}
M_{m,n}^{(1)}&=\sum_{r=1}^\infty \Omega_r^{(1)}\,
\mathcal C_{m,n}(\ba q^{-r}),\label{massagg1}\\*
\Omega_r^{(1)}&=
\frac{(\ba,1/\ba,q\ga,q\ga/\ba^2;q)_\infty}
{(q;q)_\infty^2(\ba^2,q/\ba^2;q)_\infty}
\left(\frac\ga\ba\right)^r
\frac{(q/\ba,1/\ga;q)_r}{(1/\ba,q\ga/\ba^2;q)_r},\qquad r\ge1,
\end{align}
\begin{align}
M_{m,n}^{(2)}&=\sum_{r=1}^\infty \Omega_r^{(2)}\,
\mathcal C_{m,n}(\ga^{-1}q^{-r}),\label{massagg2}\\*
\Omega_r^{(2)}&=
\frac{(1/\ga,\ga,q/\ba,q\ga^2/\ba;q)_\infty}
{(q;q)_\infty^2(1/\ga^2,q\ga^2;q)_\infty}
\left(\frac\ga\ba\right)^r
\frac{(q\ga,\ba;q)_r}{(\ga,q\ga^2/\ba;q)_r},\qquad r\ge1.
\end{align}
\end{subequations}

Only two technical points remain before the orthogonality theorem can be
stated: the interpretation of the residue sums outside their literal
convergence domains, and the check that, on the terminating lattice, these
residues are exactly the ordinary mass points.  We record both points briefly;
they are included so that the theorem has a precise meromorphic meaning.

We shall use the following finite-part convention.  If a tail has the form
\begin{equation}\label{qgeotail}
 t_r=\Lambda^r F(q^r),\qquad r\ge R,
\end{equation}
where $F(u)=\sum_{\ell\ge0}a_{\ell}u^{\ell}$ is holomorphic near $u=0$
and $R$ is sufficiently large, then, away from the polar hyperplanes
$\Lambda q^{\ell}=1$, we put
\begin{equation}\label{fptaildef}
 \operatorname{FP}\sum_{r=R}^{\infty}t_r
 :=\sum_{\ell\ge0}
 a_{\ell}\frac{(\Lambda q^{\ell})^R}{1-\Lambda q^{\ell}},
\end{equation}
where the right-hand side converges normally, and elsewhere by meromorphic
continuation.  For a full series one adds the ordinary initial segment
$\sum_{r=1}^{R-1}t_r$.  This finite part is independent of the chosen
sufficiently large $R$ and agrees with the ordinary sum when $|\Lambda|<1$.

For fixed $m,n\in\Z$, this convention applies directly to the residue
aggregates.  If $m+n$ is odd, then
$\mathcal C_{m,n}(a)=0$, hence both mass aggregates vanish.  If $m+n$ is even
and $\sigma=(m+n)/2$, then, for $j=1,2$, the $r$th residue terms in
\eqref{massaggs} have, for all sufficiently large $r$, the form
\begin{equation}\label{residue-tail-main}
 \Omega_r^{(j)}\mathcal C_{m,n}(a_r^{(j)})
 =\Lambda_j^rF_j(q^r),
\end{equation}
where $F_j$ is holomorphic at the origin and
\begin{equation}\label{mass-tail-lambdas-main}
 a_r^{(1)}=\ba q^{-r},\qquad a_r^{(2)}=\ga^{-1}q^{-r},\qquad
 \Lambda_1=\ba\ga q^\sigma,\qquad
 \Lambda_2=\frac{q^{-\sigma}}{\ba\ga}.
\end{equation}
Indeed, the parity assertion follows from
$C_k(-x;\ba,\ga\,|\,q)=(-1)^kC_k(x;\ba,\ga\,|\,q)$.  For even $m+n$, take
$a=\ba q^{-r}$ or $a=\ga^{-1}q^{-r}$, write
$x_a=(z+z^{-1})/2$ with $z=a^{1/2}$, and separate the large-$r$ dependence by
using $(u;q)_r=(u;q)_\infty/(uq^r;q)_\infty$ and
$(uq^{-r};q)_N=(-u)^Nq^{-rN+N(N-1)/2}(q^{r+1}/u;q)_N$.  The values of $C_m$
and $C_n$ on the mass lattice are handled by the convergent two-term
continuations from Section~\ref{sec:bilf}; after the finite powers of $z$ are
collected, the remaining dependence on $r$ is holomorphic in $q^r$ at the
origin.  Consequently $M_{m,n}^{(1)}$ and $M_{m,n}^{(2)}$ converge locally
uniformly in the subregions $|\Lambda_1|<1$ and $|\Lambda_2|<1$, respectively,
and elsewhere are interpreted by the finite part above.

In the sequel $M_{m,n}^{(1)}$ and $M_{m,n}^{(2)}$ denote these finite-part
sums.  Exceptional values, where a displayed denominator vanishes, a factor
$1-\Lambda_jq^\ell$ vanishes in the finite part, or mass points coalesce, are
excluded; limiting values are understood whenever the corresponding limits
exist.

Write
\begin{equation*}
\vartheta(u;q)=(u,q/u;q)_\infty.
\end{equation*}
Define
\begin{equation}\label{Hone}
H_n^{(1)}=\frac{(q,\ba,q\ga/\ba;q)_\infty^2(q\ga^2;q)_\infty}
  {(q\ga;q)_\infty^4(q/\ba;q)_\infty^2(\ba^2;q)_\infty}
\frac{(\ba^2;q)_n}{(q\ga^2;q)_n}\frac{1}{(1-\ba\ga q^{n})}
\end{equation}
and
\begin{equation}\label{Htwo}
H_n^{(2)}=\left(\frac{\ba}{q\ga}\right)^{2n}
H_{-n}^{(1)}\big|_{(\ba,\ga)\mapsto(1/\ga,1/\ba)}.
\end{equation}
Using
\begin{equation*}
(a;q)_{-n}=\frac{(-a)^{-n}q^{n(n+1)/2}}{(q/a;q)_n},
\end{equation*}
this dual definition can be written in the same form as \eqref{Hone}:
\begin{align}
H_n^{(2)}
&=-\ba\ga\,
\frac{(q,1/\ga,q\ga/\ba;q)_\infty^2(q/\ba^2;q)_\infty}
  {(q/\ba;q)_\infty^4(q\ga;q)_\infty^2(1/\ga^2;q)_\infty}
\frac{(\ba^2;q)_n}{(q\ga^2;q)_n}
\frac{1}{(1-\ba\ga q^{n})}\notag\\*
&=\ba\ga\,
\frac{\vartheta(\ga;q)^2\vartheta(\ba^2;q)}
{\vartheta(\ba;q)^2\vartheta(\ga^2;q)}\,H_n^{(1)}.
\label{Htwo-ratio}
\end{align}
In particular, the quotient \(H_n^{(2)}/H_n^{(1)}\) is independent of \(n\);
the last equality uses \(\vartheta(1/u;q)=-u^{-1}\vartheta(u;q)\).

\begin{theorem}\label{thm:orth}
Let $m,n\in\Z$.  Assume that the parameters are non-exceptional: no denominator
below vanishes and no mass points coalesce.
For $|q\ga/\ba|<1$, the first two-factor orthogonality relation is
\begin{align}
&\frac 1{2\pi}\int_{-1}^1C_m(x;\ba,\ga\,|\,q)\,C_n(x;\ba,\ga\,|\,q)
  \frac{(e^{\pm 2\ri\ta},q\ga e^{\pm 2\ri\ta}/\ba;q)_\infty}
  {(\ba e^{\pm 2\ri\ta},qe^{\pm 2\ri\ta}/\ba;q)_\infty}
\frac{\rd x}{\sqrt{1-x^2}}+M_{m,n}^{(1)}\notag\\*
&\qquad=H_n^{(1)}\,\da_{m,n}.
\label{orth-id}
\end{align}
In the same range $|q\ga/\ba|<1$, one likewise has
\begin{align}
&\frac 1{2\pi}\int_{-1}^1C_m(x;\ba,\ga\,|\,q)\,C_n(x;\ba,\ga\,|\,q)
  \frac{(e^{\pm 2\ri\ta},q\ga e^{\pm 2\ri\ta}/\ba;q)_\infty}
  {(e^{\pm 2\ri\ta}/\ga,q\ga e^{\pm 2\ri\ta};q)_\infty}
\frac{\rd x}{\sqrt{1-x^2}}+M_{m,n}^{(2)}\notag\\*
&\qquad=H_n^{(2)}\,\da_{m,n}.
\label{orth-id-neg}
\end{align}
Here the mass aggregates are interpreted as above, and exceptional parameter
values are reached, when possible, by taking limits.
\end{theorem}

\begin{proof}
We give the details for the first relation.  Start in the real terminating
subdomain $0<q<1$, $\ba>0$, $\ba q^{-s}$ non-exceptional, and $\ga=q^s$ with
$s\in\Z_{\ge0}$.  Put $\alpha=\ba q^{-s}$.  On this lattice the factor
$(1/\ga;q)_r=(q^{-s};q)_r$ makes $M_{m,n}^{(1)}$ truncate after
$r=s$, and the mass locations agree because
$\alpha q^{s-r}=\ba q^{-r}$.  With $c_s=\ba^{2s}q^{-s(s+1)}$ and
$A_s=((q;q)_s/(\alpha;q)_s)^2
=q^{s(s+1)}(q;q)_s^2/(\ba^{2s}(q/\ba;q)_s^2)$, the reduction
\begin{equation*}
 C_k(x;\ba,q^s\,|\,q)=A_sC_{k+2s}(x;\alpha\,|\,q)
\end{equation*}
and the weight identity
\begin{equation*}
 \frac{(e^{\pm2\ri\ta},q^{s+1}e^{\pm2\ri\ta}/\ba;q)_\infty}
 { (\ba e^{\pm2\ri\ta},qe^{\pm2\ri\ta}/\ba;q)_\infty}
 =c_s\frac{(e^{\pm2\ri\ta};q)_\infty}
 { (\alpha e^{\pm2\ri\ta};q)_\infty}
\end{equation*}
convert the continuous part of \eqref{orth-id} into $c_sA_s^2$ times the
continuous part of the ordinary orthogonality relation \eqref{orthmass} with
parameter $\alpha$.  Since \eqref{Cmndef} contains one half of the paired
value whereas \eqref{orthmass} uses the whole pair, direct substitution in
\eqref{massaggs} and \eqref{masspoints-lambda} gives, for $r=1,\ldots,s$ and
$j=s-r$,
\begin{equation}\label{Omega-lambda-match}
 \Omega_r^{(1)}(\ba,q^s)=2c_s\,\lambda_j(\alpha),
\end{equation}
where $\lambda_j$ is the ordinary mass weight with $\ba$ replaced by
$\alpha$.  The same simplification in \eqref{Hone} gives
\begin{equation}\label{H-terminating-match}
 H_n^{(1)}(\ba,q^s)=c_sA_s^2\,h_{n+2s}(\alpha),
\end{equation}
where $h_N(\alpha)$ is the norm in \eqref{orthmass}.
If $\alpha>1$, the ordinary relation also contains its finite mass points.
By the normalization just stated, these mass points are exactly the terms
$r=1,\ldots,s$ of $M_{m,n}^{(1)}$; if $\alpha$ lies in the mass-free range,
both sides of this statement have the empty mass sum.  Hence
\eqref{orth-id} holds for every terminating value $\ga=q^s$ in this real
subdomain, with right-hand side $c_sA_s^2h_{n+2s}(\alpha)\delta_{m,n}$, which
is \eqref{Hone} by \eqref{H-terminating-match}.

It remains to continue in $\ga$.  Fix $\ba$ and the indices $m,n$, and remove
from a small disk about $\ga=0$ the exceptional analytic hypersurfaces described
above.  In this puncture-free neighborhood the integral is analytic in $\ga$ by
uniform convergence of the products and of the bilateral series on compact
subsets whose contours stay away from poles.  The norm \eqref{Hone} is
meromorphic.  The finite-part mass aggregate is meromorphic by
the finite-part convention and the tail form \eqref{residue-tail-main}.
The apparent singularity at $\ga=0$ in the factors $(1/\ga;q)_r$ is
removable term by term after multiplication by $(\ga/\ba)^r$, since
\begin{equation*}
 \left(\frac{\ga}{\ba}\right)^r(1/\ga;q)_r
 =\ba^{-r}\prod_{\nu=0}^{r-1}(\ga-q^\nu).
\end{equation*}
Thus the difference between the two sides of \eqref{orth-id} is meromorphic
near $\ga=0$ and analytic after multiplication by a product of the finitely
many local denominator factors that occur for the fixed indices.  It vanishes
for the sequence $\ga=q^s$, $s=0,1,2,\ldots$, which has the accumulation point
$0$.  The identity theorem therefore gives \eqref{orth-id} in that local
component.  Meromorphic continuation along paths avoiding the exceptional set
then gives the stated non-exceptional parameter range $|q\ga/\ba|<1$.

For the second relation, apply the first one to the transformed parameters
$(\ba,\ga)\mapsto(1/\ga,1/\ba)$ and use the symmetry \eqref{sym}.  The
range $|q\ga/\ba|<1$ is unchanged under this symmetry.
The transformed first mass lattice is precisely the lattice in
\eqref{massagg2}, and the norm becomes \eqref{Htwo}.  This proves
\eqref{orth-id-neg}.
\end{proof}

\subsection{Four-factor orthogonality and positivity}

Set
\begin{align*}
K_{\ba}&=\frac1\ba\,
\frac{\vartheta(1/\ga;q)\vartheta(\ba^2\ga;q)}
{\vartheta(\ba/\ga;q)\vartheta(\ba\ga;q)},\\*
K_{\ga}&=\frac1{\ba^2\ga}\,
\frac{\vartheta(1/\ba\ga^2;q)\vartheta(\ba;q)}
{\vartheta(\ba/\ga;q)\vartheta(1/\ba\ga;q)},
\end{align*}
and define
\begin{equation}\label{massagg4}
M_{m,n}^{(4)}=K_{\ba}M_{m,n}^{(1)}+K_{\ga}M_{m,n}^{(2)}.
\end{equation}
Finally, put
\begin{equation}\label{Hfour}
H_n^{(4)}=-\frac{\vartheta(\ga;q)^2\vartheta(\ba^2\ga^2;q)}
{\ba\,\vartheta(\ba\ga;q)^2\vartheta(\ga^2;q)}\,H_n^{(1)}.
\end{equation}

\begin{corollary}[Four-factor orthogonality]\label{cor:fourfactor}
Assume
\begin{equation*}
|q\ga/\ba|<1.
\end{equation*}
Then, for $m,n\in\Z$, the bilateral $q$-ultraspherical functions satisfy the
full four-factor orthogonality relation
\begin{align}
&\frac 1{2\pi}\int_{-1}^1
C_m(x;\ba,\ga\,|\,q)C_n(x;\ba,\ga\,|\,q)\notag\\*
&\quad\times
  \frac{(e^{\pm 2\ri\ta},q\ga e^{\pm 2\ri\ta}/\ba,
  e^{\pm 2\ri\ta}/\ba\ga,q\ba\ga e^{\pm 2\ri\ta};q)_\infty}
  {(\ba e^{\pm 2\ri\ta},qe^{\pm 2\ri\ta}/\ba,
  e^{\pm 2\ri\ta}/\ga,q\ga e^{\pm 2\ri\ta};q)_\infty}
\frac{\rd x}{\sqrt{1-x^2}}+M_{m,n}^{(4)}
=H_n^{(4)}\,\da_{m,n}.
\label{orth-id-four}
\end{align}
\end{corollary}
\begin{proof}
Write $y=e^{2\ri\ta}$ and denote by $W_1(y)$, $W_2(y)$, and $W_4(y)$ the
continuous weights in \eqref{orth-id}, \eqref{orth-id-neg}, and
\eqref{orth-id-four}, respectively.  Thus
\begin{align*}
W_1(y)&=\frac{(y^{\pm1},q\ga y^{\pm1}/\ba;q)_\infty}
{(\ba y^{\pm1},qy^{\pm1}/\ba;q)_\infty},\qquad
W_2(y)=\frac{(y^{\pm1},q\ga y^{\pm1}/\ba;q)_\infty}
{(y^{\pm1}/\ga,q\ga y^{\pm1};q)_\infty},\\*
W_4(y)&=\frac{(y^{\pm1},q\ga y^{\pm1}/\ba,
y^{\pm1}/\ba\ga,q\ba\ga y^{\pm1};q)_\infty}
{(\ba y^{\pm1},qy^{\pm1}/\ba,
y^{\pm1}/\ga,q\ga y^{\pm1};q)_\infty}.
\end{align*}
The Weierstrass addition formula for theta functions, in the form
\begin{align}
&K_{\ba}\,\vartheta(y/\ga;q)\vartheta(q\ga y;q)
+K_{\ga}\,\vartheta(\ba y;q)\vartheta(qy/\ba;q)\notag\\*
&\qquad=\vartheta(y/\ba\ga;q)\vartheta(q\ba\ga y;q),
\label{theta-four}
\end{align}
implies
\begin{equation}\label{weight-four}
W_4(y)=K_{\ba}W_1(y)+K_{\ga}W_2(y).
\end{equation}
Indeed, after the common factor
\begin{equation*}
\frac{(y^{\pm1},q\ga y^{\pm1}/\ba;q)_\infty}
{(\ba y^{\pm1},qy^{\pm1}/\ba,
y^{\pm1}/\ga,q\ga y^{\pm1};q)_\infty}
\end{equation*}
is factored out of the right-hand side of \eqref{weight-four}, the remaining
identity is precisely \eqref{theta-four}, since
\begin{align*}
(y^{\pm1}/\ga,q\ga y^{\pm1};q)_\infty
&=\vartheta(y/\ga;q)\vartheta(q\ga y;q),\\*
(\ba y^{\pm1},qy^{\pm1}/\ba;q)_\infty
&=\vartheta(\ba y;q)\vartheta(qy/\ba;q),\\*
(y^{\pm1}/\ba\ga,q\ba\ga y^{\pm1};q)_\infty
&=\vartheta(y/\ba\ga;q)\vartheta(q\ba\ga y;q).
\end{align*}

Let $I_{m,n}^{(j)}$ be the integral in \eqref{orth-id},
\eqref{orth-id-neg}, or \eqref{orth-id-four} with weight $W_j$, for
$j=1,2,4$.  Multiplying \eqref{weight-four} by
$C_m(x;\ba,\ga\,|\,q)C_n(x;\ba,\ga\,|\,q)$ and integrating gives
\begin{equation*}
I_{m,n}^{(4)}=K_{\ba}I_{m,n}^{(1)}+K_{\ga}I_{m,n}^{(2)}.
\end{equation*}
Together with \eqref{massagg4}, \eqref{orth-id}, and \eqref{orth-id-neg}, this
gives
\begin{align*}
I_{m,n}^{(4)}+M_{m,n}^{(4)}
&=K_{\ba}\bigl(I_{m,n}^{(1)}+M_{m,n}^{(1)}\bigr)
 +K_{\ga}\bigl(I_{m,n}^{(2)}+M_{m,n}^{(2)}\bigr)\\*
&=\bigl(K_{\ba}H_n^{(1)}+K_{\ga}H_n^{(2)}\bigr)\da_{m,n}.
\end{align*}
It remains only to simplify the norm.  Substitution of \eqref{Htwo-ratio}
and the definitions of $K_{\ba}$ and
$K_{\ga}$ gives
\begin{equation*}
K_{\ba}H_n^{(1)}+K_{\ga}H_n^{(2)}
=-\frac{\vartheta(\ga;q)^2\vartheta(\ba^2\ga^2;q)}
{\ba\,\vartheta(\ba\ga;q)^2\vartheta(\ga^2;q)}\,H_n^{(1)},
\end{equation*}
where the last equality is again the Weierstrass addition formula, now in
the corresponding specialization of the parameters.  By \eqref{Hfour}, this is
$H_n^{(4)}$, proving \eqref{orth-id-four}.
\end{proof}

\begin{remark}[Positivity]\label{rem:bilateral-positivity}
For real parameters and \(0<q<1\), the continuous factors \(W_1\), \(W_2\),
and \(W_4\) are pointwise non-negative on the unit circle whenever their
parameters are real and no denominator factor has a zero on the contour
(equivalently, each real denominator parameter \(a\) in a factor
\((ay^{\pm1};q)_\infty\) avoids \(\pm q^{-\Z_{\ge0}}\)).  This pointwise
condition is not, however, the same as positive definiteness of the full
orthogonality functional, because the residue aggregates must also be included.

A necessary Favard-type condition, and the natural bilateral analogue of the
Askey--Ismail condition displayed in Section~\ref{sec:c}, is obtained directly
from the three-term recurrence \eqref{rec}.  Any non-degenerate positive
orthogonality functional for the bilateral sequence must satisfy
\begin{equation}\label{bil-pos-cond}
\frac{(1-\ba^2q^n)(1-\ga^2q^{n+1})}
{(1-\ba\ga q^n)(1-\ba\ga q^{n+1})}>0,
\qquad n\in\Z,
\end{equation}
and the initial norm \(H_0^{(j)}\), for the particular normalization of the
functional under consideration, must have positive sign.  Conversely,
\eqref{bil-pos-cond}, together with a positive initial norm, is the formal
recurrence-level positivity condition after symmetrizing the corresponding
doubly infinite Jacobi matrix.
In the specialization \(\ga=1\), where the negative-order part collapses and
one restricts to the ordinary Rogers index set, \eqref{bil-pos-cond} reduces to
the condition quoted in Section~\ref{sec:c}.  For genuinely bilateral
parameters the inequalities in \eqref{bil-pos-cond} depend on the position of
\(\ba\), \(\ga\), and \(\ba\ga\) relative to the two-sided \(q\)-lattice, and
there does not seem to be a single interval description comparable to the
Rogers case.  Thus we use \eqref{bil-pos-cond} as the usable full condition
at the recurrence level; a simple measure-theoretic description of the
positive parameter domains for the explicit integral-plus-residue realizations
remains more delicate.
\end{remark}

\subsection{Integral evaluations deduced from the orthogonality relations}

The preceding orthogonality relations also imply integral evaluations with
generated mass aggregates.  These are obtained by applying the bilateral
generating functions from Theorems~\ref{thmgf} and \ref{thmgf2}.  For
\(\epsilon=0,1\), put
\begin{equation}\label{genkernels}
\mathcal G_\epsilon(x;t)=
\sum_{n\in\Z}(1-\ba\ga q^n)^\epsilon
C_n(x;\ba,\ga\,|\,q)t^n.
\end{equation}
Thus \(\mathcal G_0\) is the generating function in \eqref{gfbil}, while
\(\mathcal G_1\) is the generating function in \eqref{gfbil2}.  For later
reference we also name the three continuous weights
\begin{subequations}\label{contweights}
\begin{align}
W_1(y)&=\frac{(y^{\pm1},q\ga y^{\pm1}/\ba;q)_\infty}
{(\ba y^{\pm1},qy^{\pm1}/\ba;q)_\infty},\label{contweight-one}\\*
W_2(y)&=\frac{(y^{\pm1},q\ga y^{\pm1}/\ba;q)_\infty}
{(y^{\pm1}/\ga,q\ga y^{\pm1};q)_\infty},\label{contweight-two}\\*
W_4(y)&=\frac{(y^{\pm1},q\ga y^{\pm1}/\ba,
 y^{\pm1}/\ba\ga,q\ba\ga y^{\pm1};q)_\infty}
{(\ba y^{\pm1},qy^{\pm1}/\ba,
 y^{\pm1}/\ga,q\ga y^{\pm1};q)_\infty}.
\label{contweight-four}
\end{align}
\end{subequations}
For \(a\ne0\), define the generating-function parity-symmetrized mass kernel
\begin{equation}\label{genmasskernel}
\mathcal E_{\epsilon,\eta}(a;t_1,t_2)=\frac12\Big(
\mathcal G_\epsilon(x_a;t_1)\mathcal G_\eta(x_a;t_2)
+\mathcal G_\epsilon(-x_a;t_1)\mathcal G_\eta(-x_a;t_2)\Big),
\end{equation}
where \(x_a=(a^{1/2}+a^{-1/2})/2\).  The corresponding generating-function
mass aggregates are
\begin{subequations}\label{genmasses}
\begin{align}
\mathcal M_1^{\epsilon,\eta}(t_1,t_2)&=
\sum_{r=1}^\infty \Omega_r^{(1)}
\mathcal E_{\epsilon,\eta}(\ba q^{-r};t_1,t_2),\\*
\mathcal M_2^{\epsilon,\eta}(t_1,t_2)&=
\sum_{r=1}^\infty \Omega_r^{(2)}
\mathcal E_{\epsilon,\eta}(\ga^{-1}q^{-r};t_1,t_2),\\*
\mathcal M_4^{\epsilon,\eta}(t_1,t_2)&=
K_{\ba}\mathcal M_1^{\epsilon,\eta}(t_1,t_2)
+K_{\ga}\mathcal M_2^{\epsilon,\eta}(t_1,t_2).
\end{align}
\end{subequations}
Finally, with \(u=t_1t_2\), let
\begin{equation}\label{gennorms}
\mathcal H_j^{[\ell]}(u)=
\sum_{n\in\Z}(1-\ba\ga q^n)^\ell H_n^{(j)}u^n,
\qquad j=1,2,4,\quad \ell=0,1,2,
\end{equation}
where \(H_n^{(4)}\) is defined in \eqref{Hfour}.  Equivalently,
\begin{align*}
\mathcal H_j^{[1]}(u)&=\mathcal H_j^{[0]}(u)-\ba\ga\mathcal H_j^{[0]}(qu),\\*
\mathcal H_j^{[2]}(u)&=\mathcal H_j^{[0]}(u)-2\ba\ga\mathcal H_j^{[0]}(qu)
+\ba^2\ga^2\mathcal H_j^{[0]}(q^2u).
\end{align*}
For compactness in the following product forms, set
\begin{equation*}
\mathcal P_{\epsilon,\eta}(t_1,t_2)=
(1-\ga/\ba)^{\epsilon+\eta}
(1-\ga/\ba t_1^2)^\epsilon
(1-\ga/\ba t_2^2)^\eta
\frac{(q,q\ga/\ba;q)_\infty^4}{(q\ga,q/\ba;q)_\infty^4}.
\end{equation*}

\begin{theorem}[Mass-aggregate integral evaluations]
\label{thm:mass-integral-evals}
Let \(\epsilon,\eta\in\{0,1\}\).  In any common domain of absolute convergence,
and elsewhere by meromorphic continuation of both sides, the following integral
evaluations, deduced from the orthogonality relations, hold under the single
parameter condition \(|q\ga/\ba|<1\).  For the first two-factor weight
\(W_1\) in \eqref{contweight-one},
\begin{align}
&\mathcal P_{\epsilon,\eta}(t_1,t_2)\frac1{2\pi}\int_{-1}^1
\frac{(\ba t_1e^{\pm\ri\ta},qe^{\pm\ri\ta}/\ba t_1,
\ba t_2e^{\pm\ri\ta},qe^{\pm\ri\ta}/\ba t_2;q)_\infty}
{(t_1e^{\pm\ri\ta},q^{1-\epsilon}\ga e^{\pm\ri\ta}/\ba t_1,
 t_2e^{\pm\ri\ta},q^{1-\eta}\ga e^{\pm\ri\ta}/\ba t_2;q)_\infty}
 W_1(e^{2\ri\ta})\frac{\rd x}{\sqrt{1-x^2}}\notag\\*
&\qquad\qquad+\mathcal M_1^{\epsilon,\eta}(t_1,t_2)
=\mathcal H_1^{[\epsilon+\eta]}(t_1t_2).
\label{genorth1}
\end{align}
For the second two-factor weight \(W_2\) in \eqref{contweight-two},
\begin{align}
&\mathcal P_{\epsilon,\eta}(t_1,t_2)\frac1{2\pi}\int_{-1}^1
\frac{(\ba t_1e^{\pm\ri\ta},qe^{\pm\ri\ta}/\ba t_1,
\ba t_2e^{\pm\ri\ta},qe^{\pm\ri\ta}/\ba t_2;q)_\infty}
{(t_1e^{\pm\ri\ta},q^{1-\epsilon}\ga e^{\pm\ri\ta}/\ba t_1,
 t_2e^{\pm\ri\ta},q^{1-\eta}\ga e^{\pm\ri\ta}/\ba t_2;q)_\infty}
 W_2(e^{2\ri\ta})\frac{\rd x}{\sqrt{1-x^2}}\notag\\*
&\qquad\qquad+\mathcal M_2^{\epsilon,\eta}(t_1,t_2)
=\mathcal H_2^{[\epsilon+\eta]}(t_1t_2).
\label{genorth2}
\end{align}
For the four-factor weight \(W_4\) in \eqref{contweight-four},
\begin{align}
&\mathcal P_{\epsilon,\eta}(t_1,t_2)\frac1{2\pi}\int_{-1}^1
\frac{(\ba t_1e^{\pm\ri\ta},qe^{\pm\ri\ta}/\ba t_1,
\ba t_2e^{\pm\ri\ta},qe^{\pm\ri\ta}/\ba t_2;q)_\infty}
{(t_1e^{\pm\ri\ta},q^{1-\epsilon}\ga e^{\pm\ri\ta}/\ba t_1,
 t_2e^{\pm\ri\ta},q^{1-\eta}\ga e^{\pm\ri\ta}/\ba t_2;q)_\infty}
 W_4(e^{2\ri\ta})\frac{\rd x}{\sqrt{1-x^2}}\notag\\*
&\qquad\qquad+\mathcal M_4^{\epsilon,\eta}(t_1,t_2)
=\mathcal H_4^{[\epsilon+\eta]}(t_1t_2).
\label{genorth4}
\end{align}
In particular, in \eqref{genorth1} the choices
\((\epsilon,\eta)=(0,0),(0,1),(1,1)\) give, respectively, the
mass-aggregate extensions of \eqref{thm1322}, \eqref{orth-idgf2-cont}, and
\eqref{orth-idgf3-cont}; the remaining mixed choice \((1,0)\) is the
\(t_1,t_2\)-interchanged version of \((0,1)\).
\end{theorem}
\begin{proof}
Multiply \eqref{orth-id}, \eqref{orth-id-neg}, or \eqref{orth-id-four} by
\((1-\ba\ga q^m)^\epsilon(1-\ba\ga q^n)^\eta t_1^m t_2^n\), sum over
\(m,n\in\Z\), and use \eqref{genkernels}.  Substitution of the product forms
in Theorems~\ref{thmgf} and \ref{thmgf2}, together with the displayed
definitions of the weights \(W_1,W_2,W_4\) in \eqref{contweights}, gives
the three product integrands above.  The parity-symmetrized products in
\eqref{Cmndef} generate precisely the kernels \eqref{genmasskernel}, giving
\eqref{genmasses}.  The right-hand side is \eqref{gennorms}.  The initial
interchanges of summation, integration, and residue summation are justified in
an absolute-convergence subdomain; the stated identities then extend
meromorphically in the parameters and variables.
\end{proof}

\subsection{Specializations and closed norm sums}

The generic generating-function mass aggregates in \eqref{genmasses} do not
seem to admit a simple product evaluation.  Several degenerate and terminating
specializations, however, are explicit.

\begin{corollary}[Degenerate generating-function mass aggregates]\label{cor:vanishing-genmasses}
Let \(\epsilon,\eta\in\{0,1\}\).  Suppose that either
\(\epsilon=1\) and \(t_1^2=\ga/\ba\), or \(\eta=1\) and
\(t_2^2=\ga/\ba\).  At non-exceptional parameter values for which the
corresponding continuous product integral is finite, the continuous term in
Theorem~\ref{thm:mass-integral-evals} is annihilated by the factor
\(\mathcal P_{\epsilon,\eta}(t_1,t_2)\).  Consequently the finite-part
generating-function mass aggregate is isolated:
\begin{align}
\mathcal M_1^{\epsilon,\eta}(t_1,t_2)&=\mathcal H_1^{[\epsilon+\eta]}(t_1t_2),
\label{genmass-degen1}\\*
\mathcal M_2^{\epsilon,\eta}(t_1,t_2)&=\mathcal H_2^{[\epsilon+\eta]}(t_1t_2),
\label{genmass-degen2}\\*
\mathcal M_4^{\epsilon,\eta}(t_1,t_2)&=\mathcal H_4^{[\epsilon+\eta]}(t_1t_2),
\label{genmass-degen4}
\end{align}
in the common range \(|q\ga/\ba|<1\).

In particular, if \(\rho^2=\ga/\ba\), then in this same range
\begin{align}
\mathcal M_j^{1,0}(\rho,t)&=\mathcal H_j^{[1]}(\rho t),
\label{genorth-vanish10}\\*
\mathcal M_j^{0,1}(t,\rho)&=\mathcal H_j^{[1]}(t\rho),
\label{genorth-vanish01}\\*
\mathcal M_j^{1,1}(\rho,t)&=\mathcal H_j^{[2]}(\rho t),
\label{genorth-vanish11a}\\*
\mathcal M_j^{1,1}(t,\rho)&=\mathcal H_j^{[2]}(t\rho),
\label{genorth-vanish11b}
\end{align}
for \(j=1,2,4\).  The specialization \(t_1=t_2=\rho\) is obtained from
either of the last two formulae by taking \(t=\rho\).  In the limiting case
\(\ga=\ba\), the corresponding finite-part statements for \(j=1,2\) are
obtained by taking the non-exceptional limit \(\rho\to1\).  More generally,
for \(j=1,2\) and \(\epsilon+\eta>0\), whenever the limit exists and the
continuous product term remains bounded,
\begin{equation}
\lim_{\ga\to\ba}\mathcal M_j^{\epsilon,\eta}(t_1,t_2)
=\lim_{\ga\to\ba}\mathcal H_j^{[\epsilon+\eta]}(t_1t_2).
\label{genorth-vanish-betagamma}
\end{equation}
The mass term in these identities is essential; in general the corresponding
mass-free integral evaluations are not valid.
\end{corollary}
\begin{proof}
Under either displayed specialization the factor
\(\mathcal P_{\epsilon,\eta}(t_1,t_2)\) vanishes.  If the continuous product
integral has a finite value at the specialization, the continuous term in
\eqref{genorth1}, \eqref{genorth2}, or \eqref{genorth4} therefore contributes
zero, and the three identities \eqref{genmass-degen1}--\eqref{genmass-degen4}
follow directly from Theorem~\ref{thm:mass-integral-evals}.  The formulae with
\(\rho^2=\ga/\ba\) are the corresponding one-variable specializations.

It is important that the product formula \eqref{gfbil2} cannot be used
pointwise at the mass points to conclude that the generating-function mass kernels in
\eqref{genmasskernel} vanish.  At such points the same specialization may
also produce zeros in denominator factors, so the finite-part value of the
mass aggregate can contain a non-zero limiting contribution.  The limiting
statement for \(\ga\to\ba\) is obtained in the same way from \eqref{genorth1}
and \eqref{genorth2}, whenever the indicated non-exceptional limit exists.
\end{proof}

\begin{remark}[Explicit residue summations behind \eqref{genmass-degen1}--\eqref{genorth-vanish11b}]
The identities \eqref{genmass-degen1}--\eqref{genorth-vanish11b} can be written
as standalone basic-hypergeometric residue summations.  Put
$Q_{\ba,\ga}=(q,q\ga/\ba;q)_\infty^2/(q\ga,q/\ba;q)_\infty^2$ and, for \(z\ne0\), define the two product kernels
\begin{align*}
\Phi_0(z;t)&=Q_{\ba,\ga}
\frac{(\ba t z^{\pm1},qz^{\pm1}/(\ba t);q)_\infty}
{(t z^{\pm1},q\ga z^{\pm1}/(\ba t);q)_\infty},\\*
\Phi_1(z;t)&=(1-\ga/\ba)(1-\ga/\ba t^2)\,Q_{\ba,\ga}
\frac{(\ba t z^{\pm1},qz^{\pm1}/(\ba t);q)_\infty}
{(t z^{\pm1},\ga z^{\pm1}/(\ba t);q)_\infty}.
\end{align*}
Thus \(\Phi_\epsilon(z;t)=\mathcal G_\epsilon((z+z^{-1})/2;t)\), by
Theorems~\ref{thmgf} and \ref{thmgf2}.  For \(a\ne0\), let
\[
\widetilde{\mathcal E}_{\epsilon,\eta}(a;t_1,t_2)=\frac12\Big(
\Phi_\epsilon(a^{1/2};t_1)\Phi_\eta(a^{1/2};t_2)
+\Phi_\epsilon(-a^{1/2};t_1)\Phi_\eta(-a^{1/2};t_2)\Big),
\]
where either branch of \(a^{1/2}\) may be used.  Finally set
\begin{align*}
\mathcal R_1^{\epsilon,\eta}(t_1,t_2)&=
\sum_{r=1}^\infty \Omega_r^{(1)}
\widetilde{\mathcal E}_{\epsilon,\eta}(\ba q^{-r};t_1,t_2),\\*
\mathcal R_2^{\epsilon,\eta}(t_1,t_2)&=
\sum_{r=1}^\infty \Omega_r^{(2)}
\widetilde{\mathcal E}_{\epsilon,\eta}(\ga^{-1}q^{-r};t_1,t_2),\\*
\mathcal R_4^{\epsilon,\eta}(t_1,t_2)&=
K_{\ba}\mathcal R_1^{\epsilon,\eta}(t_1,t_2)
+K_{\ga}\mathcal R_2^{\epsilon,\eta}(t_1,t_2).
\end{align*}
These are precisely the mass aggregates \(\mathcal M_j^{\epsilon,\eta}\), but
with the generating functions evaluated as explicit products.  If the products
\(\Phi_\epsilon\) and \(\Phi_\eta\) are replaced by their defining bilateral
series, the same left-hand sides become double bilateral sums with the
additional residue summation over \(r\).

The product right-hand sides can be written uniformly as follows.  Define
\begin{align*}
P_1(u)&=
\frac{(q,\ba,q\ga/\ba;q)_\infty^2(q\ga^2;q)_\infty}
{(q\ga;q)_\infty^4(q/\ba;q)_\infty^2(\ba^2;q)_\infty}
\frac{(q,\ba^2u,q/\ba^2u,q\ga^2/\ba^2;q)_\infty}
{(q\ga^2,u,q\ga^2/\ba^2u,q/\ba^2;q)_\infty},\\*
P_2(u)&=P_1(u)\frac{(\ga-\ba)(\ga+\ba u)}{(\ga^2-\ba^2u)},
\end{align*}
and constants
\[
D_1=1,\qquad
D_2=\ba\ga\,
\frac{\vartheta(\ga;q)^2\vartheta(\ba^2;q)}
{\vartheta(\ba;q)^2\vartheta(\ga^2;q)},\qquad
D_4=-\frac{\vartheta(\ga;q)^2\vartheta(\ba^2\ga^2;q)}
{\ba\,\vartheta(\ba\ga;q)^2\vartheta(\ga^2;q)}.
\]
Then the explicit version of \eqref{genmass-degen1}--\eqref{genmass-degen4} is
$\mathcal R_j^{\epsilon,\eta}(t_1,t_2)=D_j P_{\epsilon+\eta}(t_1t_2)$ for $j=1,2,4$, where \(\epsilon+\eta\in\{1,2\}\), whenever either \(\epsilon=1\) and
\(t_1^2=\ga/\ba\), or \(\eta=1\) and
\(t_2^2=\ga/\ba\), with the same finite-part interpretation as in
Corollary~\ref{cor:vanishing-genmasses}.  In particular, if \(\rho^2=\ga/\ba\),
then \eqref{genorth-vanish10}--\eqref{genorth-vanish11b} become the four
product identities
\begin{align*}
\mathcal R_j^{1,0}(\rho,t)&=D_jP_1(\rho t),&
\mathcal R_j^{0,1}(t,\rho)&=D_jP_1(t\rho),\\*
\mathcal R_j^{1,1}(\rho,t)&=D_jP_2(\rho t),&
\mathcal R_j^{1,1}(t,\rho)&=D_jP_2(t\rho),
\qquad j=1,2,4.
\end{align*}
For \(j=1\), these product sides are Ramanujan's \({}_1\psi_1\) summation
and its first contiguous difference applied to the norm-generating functions;
for \(j=2\) and \(j=4\), they are obtained from the same products by the
factors \eqref{Htwo-ratio} and \eqref{Hfour}.  The nontrivial point, when read
purely as basic hypergeometric identities, is that the explicit residue sums
\(\mathcal R_j^{\epsilon,\eta}\) collapse to these products.
\end{remark}

\begin{corollary}[Terminating generating-function mass aggregates]\label{cor:terminating-genmasses}
Let \(s\in\Z_{\ge0}\).  If \(\ga=q^s\), then
\begin{equation*}
\mathcal M_1^{\epsilon,\eta}(t_1,t_2)=
\sum_{r=1}^s \Omega_r^{(1)}
\mathcal E_{\epsilon,\eta}(\ba q^{-r};t_1,t_2),
\end{equation*}
with the empty sum interpreted as zero when \(s=0\).  If \(\ba=q^{-s}\), then
\begin{equation*}
\mathcal M_2^{\epsilon,\eta}(t_1,t_2)=
\sum_{r=1}^s \Omega_r^{(2)}
\mathcal E_{\epsilon,\eta}(\ga^{-1}q^{-r};t_1,t_2).
\end{equation*}
Consequently \(\mathcal M_4^{\epsilon,\eta}\) is finite whenever both
component aggregates terminate and the coefficients \(K_{\ba}\) and
\(K_{\ga}\) are finite.
\end{corollary}
\begin{proof}
For \(\ga=q^s\), the factor \((1/\ga;q)_r=(q^{-s};q)_r\) in
\(\Omega_r^{(1)}\) vanishes for every \(r>s\).  Similarly, for
\(\ba=q^{-s}\), the factor \((\ba;q)_r=(q^{-s};q)_r\) in
\(\Omega_r^{(2)}\) vanishes for every \(r>s\).  The statement for
\(\mathcal M_4\) follows from \eqref{genmasses}.
\end{proof}

\begin{corollary}[Closed norm-generating functions]\label{cor:closednorms}
Set
\begin{equation}\label{normconsts}
B_{\ba,\ga}=\frac{(q,\ba,q\ga/\ba;q)_\infty^2(q\ga^2;q)_\infty}
{(q\ga;q)_\infty^4(q/\ba;q)_\infty^2(\ba^2;q)_\infty},
\qquad
A_{\ba,\ga}=\frac{B_{\ba,\ga}}{(1-\ba\ga)}.
\end{equation}
Then
\begin{align}
\mathcal H_1^{[0]}(u)&=
A_{\ba,\ga}
{}_2\psi_2\!\left[\begin{matrix}\ba^2,\ba\ga\\
q\ga^2,q\ba\ga\end{matrix};q,u\right],
\label{Hgen10}\\*
\mathcal H_1^{[1]}(u)&=
B_{\ba,\ga}
\frac{(q,\ba^2u,q/\ba^2u,q\ga^2/\ba^2;q)_\infty}
{(q\ga^2,u,q\ga^2/\ba^2u,q/\ba^2;q)_\infty},
\label{Hgen11}\\*
\mathcal H_1^{[2]}(u)&=
B_{\ba,\ga}
\frac{(q,\ba^2u,q/\ba^2u,q\ga^2/\ba^2;q)_\infty}
{(q\ga^2,u,q\ga^2/\ba^2u,q/\ba^2;q)_\infty}
\frac{(\ga-\ba)(\ga+\ba u)}{(\ga^2-\ba^2u)}.
\label{Hgen12}
\end{align}
Moreover,
\begin{align}
\mathcal H_2^{[0]}(u)&=
A_{1/\ga,1/\ba}
{}_2\psi_2\!\left[\begin{matrix}1/\ga^2,1/\ba\ga\\
q/\ba^2,q/\ba\ga\end{matrix};q,
\frac{q^2\ga^2}{\ba^2u}\right],
\label{Hgen20}\\*
\mathcal H_4^{[\ell]}(u)&=
-\frac{\vartheta(\ga;q)^2\vartheta(\ba^2\ga^2;q)}
{\ba\,\vartheta(\ba\ga;q)^2\vartheta(\ga^2;q)}
\mathcal H_1^{[\ell]}(u),\qquad \ell=0,1,2,
\label{Hgen4}
\end{align}
and \(\mathcal H_2^{[1]}\), \(\mathcal H_2^{[2]}\) are obtained from
\eqref{gennorms} by the two finite-difference formulae following it.
\end{corollary}
\begin{proof}
Formula \eqref{Hgen10} follows from \eqref{Hone} and
\((\ba\ga;q)_n/(q\ba\ga;q)_n=(1-\ba\ga)/(1-\ba\ga q^n)\).  Multiplying
\eqref{Hone} by \((1-\ba\ga q^n)\) gives a \({}_1\psi_1\) series, and
Ramanujan's summation \eqref{1psi1} gives \eqref{Hgen11}.  Multiplication by
\((1-\ba\ga q^n)^2\) gives
\begin{equation*}
\mathcal H_1^{[2]}(u)=B_{\ba,\ga}
\left(S(u)-\ba\ga S(qu)\right),
\qquad
S(u)={}_1\psi_1\!\left[\begin{matrix}\ba^2\\ q\ga^2\end{matrix};q,u\right].
\end{equation*}
Using \eqref{1psi1} and
\begin{equation*}
\frac{S(qu)}{S(u)}=\frac{(1-u)}{(\ga^2-\ba^2u)}
\end{equation*}
then yields \eqref{Hgen12}.  Formula \eqref{Hgen20} follows from
\eqref{Htwo} by setting \(k=-n\) in the generating-function norm sum.  Finally,
\eqref{Hgen4} follows directly from \eqref{Hfour}.
\end{proof}

\begin{example}[Residue identities from the integral evaluations]\label{cor:residue-degenerate-masses}\label{ex:residue-integral-evals}
Let \(\rho^2=\ga/\ba\), put \(u=\rho t\), and set
\[
S_{\ba,\ga}(u)={}_1\psi_1\!\left[
\begin{matrix}\ba^2\\ q\ga^2\end{matrix};q,u\right],
\qquad
T_{\ba,\ga}(u)={}_2\psi_2\!\left[
\begin{matrix}1/\ga^2,1/\ba\ga\\ q/\ba^2,q/\ba\ga\end{matrix};q,
\frac{q^2\ga^2}{\ba^2u}\right].
\]
Putting \(t_1=\rho\) and \(\epsilon=1\) in Theorem~\ref{thm:mass-integral-evals}
(or symmetrically \(t_2=\rho\), \(\eta=1\)) makes
\(\mathcal P_{\epsilon,\eta}(t_1,t_2)\) vanish.  Moving the contour past
the pole lattice therefore isolates the finite-part residue contribution.
For the first residue family this gives Ramanujan's summation and its first
contiguous difference:
\begin{subequations}\label{residue-degen-first}
\begin{align}
S_{\ba,\ga}(u)&=
\frac{(q,\ba^2u,q/\ba^2u,q\ga^2/\ba^2;q)_\infty}
{(q\ga^2,u,q\ga^2/\ba^2u,q/\ba^2;q)_\infty},\label{residue-degen-first1}\\*
S_{\ba,\ga}(u)-\ba\ga S_{\ba,\ga}(qu)&=
\frac{(q,\ba^2u,q/\ba^2u,q\ga^2/\ba^2;q)_\infty}
{(q\ga^2,u,q\ga^2/\ba^2u,q/\ba^2;q)_\infty}
\frac{(\ga-\ba)(\ga+\ba u)}{(\ga^2-\ba^2u)}.
\label{residue-degen-first2}
\end{align}
\end{subequations}
The second residue family gives, for instance, the contiguous \({}_2\psi_2\) transformation
\begin{align}
A_{1/\ga,1/\ba}\big(T_{\ba,\ga}(u)-\ba\ga T_{\ba,\ga}(qu)\big)
&=\Lambda_{\ba,\ga}B_{\ba,\ga}
\frac{(q,\ba^2u,q/\ba^2u,q\ga^2/\ba^2;q)_\infty}
{(q\ga^2,u,q\ga^2/\ba^2u,q/\ba^2;q)_\infty},
\label{residue-degen-second1}
\end{align}
where
$\Lambda_{\ba,\ga}=\ba\ga(\ga,q/\ga;q)_\infty^2(\ba^2,q/\ba^2;q)_\infty/
((\ba,q/\ba;q)_\infty^2(\ga^2,q/\ga^2;q)_\infty)$.  The second finite difference of \(T_{\ba,\ga}\) and the four-factor residue identity
are obtained in the same way from \(\mathcal H_2^{[2]}\) and the theta-linear
combination \eqref{Hgen4}; we do not record them separately.
\end{example}
\begin{corollary}[Classical Rogers specializations]\label{cor:rogers-integral-evals}
In the specialization \(\ga=1\), the first integral evaluation
\eqref{genorth1} with \((\epsilon,\eta)=(0,0)\) reduces to the classical
integral evaluation \eqref{thm1322}.  With \((\epsilon,\eta)=(0,1)\) and
\((1,1)\), respectively, and after using the limiting second generating
function \eqref{cgfbil2} and dividing by its explicit scalar factors
\((1-\ba)(1-\ba t_2^2)\) and
\((1-\ba)^2(1-\ba t_1^2)(1-\ba t_2^2)\), it reduces to the two evaluations
\eqref{orth-idgf2-cont} and \eqref{orth-idgf3-cont}.
\end{corollary}
\begin{proof}
For \(\ga=1\), the bilateral functions with negative index vanish and
\(C_n(x;\ba,1\,|\,q)=C_n(x;\ba\,|\,q)\) for \(n\ge0\).  Also
\((1/\ga;q)_r=(1;q)_r=0\) for every \(r\ge1\), so the generating-function
first mass aggregate disappears.  The first generating function \eqref{gfbil} becomes
\eqref{gfuni}, while \eqref{gfbil2} becomes \eqref{cgfbil2}.  Substituting
these limits in \eqref{genorth1} and simplifying the right-hand sides in
Corollary~\ref{cor:closednorms} gives exactly \eqref{thm1322},
\eqref{orth-idgf2-cont}, and \eqref{orth-idgf3-cont}.
\end{proof}

\section{Mixed and quasi-linearization formulae, and mass aggregates}
\label{sec:linear}

In this section we prove mixed and quasi-linearization formulae and record the
resulting identities for the residue mass aggregates.  In the bilateral
setting, product formulae for the functions immediately produce identities for
these aggregates.  This separates the genuinely finite situations from the
quasi-linearized ones, in which an additional analytic remainder remains
visible.

There are two natural unilateral factors.  The parameter $\ba/\ga$ is adapted
to the bilateral family and gives a finite expansion in the same basis.  The
classical Rogers parameter $\ba$ is more symmetric in the limit $\ga\to1$, but
for generic $\ga$ it gives only a quasi-linearization: an infinite Rogers-type
part plus an explicit analytic remainder.

\subsection{A finite mixed linearization formula}

\begin{theorem}[Mixed bilateral--Rogers linearization]\label{thm:mixedlin}
Let $m\in\Z$, $n\in\Z_{\ge0}$, and write
$x=(z+z^{-1})/2$.  If the bilateral functions in the identity below are
interpreted by their defining bilateral series, assume the common convergence
condition
\begin{equation*}
\left|\frac{q\ga z^2}{\ba}\right|<1,
\qquad
\left|\frac{q\ga}{\ba z^2}\right|<1,
\end{equation*}
equivalently $|q\ga/\ba|<|z|^2<|\ba/q\ga|$.
For $z$ on the unit circle this reduces to $|q\ga/\ba|<1$.
Away from exceptional parameter values at which denominators vanish, we have
\begin{equation}\label{mixedlinear}
C_m(x;\ba,\ga\,|\,q)C_n(x;\ba/\ga\,|\,q)
=\sum_{j=0}^{n}L_{m,n}^{(j)}(\ba,\ga\,|\,q)
 C_{m-n+2j}(x;\ba,\ga\,|\,q),
\end{equation}
where
\begin{align}
L_{m,n}^{(j)}(\ba,\ga\,|\,q)
&=\frac{(\ba/\ga;q)_n}{(q;q)_n}
\frac{(\ba^2q^{m-n};q)_n}{(\ba\ga q^{m-n+1};q)_n}
\frac{(1-\ba\ga q^{m-n+2j})}{(1-\ba\ga q^{m-n})}\notag\\*
&\quad\times
\frac{(\ba\ga q^{m-n},\ba/\ga,\ba^2q^m,q^{-n};q)_j}
{(q,\ga^2q^{m-n+1},\ga q^{1-n}/\ba,\ba\ga q^{m+1};q)_j}
\frac{(\ga^2q^{m-n+1};q)_{2j}}{(\ba^2q^{m-n};q)_{2j}}
\left(\frac{\ga q}{\ba}\right)^j.
\label{mixedcoeff}
\end{align}
The preceding convergence restrictions are only needed for the direct
termwise interpretation of the bilateral series.  Both sides are meromorphic
in $\ba$ and $\ga$, so the identity extends by meromorphic continuation to all
non-exceptional parameters.
\end{theorem}

\begin{proof}
Write $x=(z+z^{-1})/2$ and compare Laurent coefficients in $z$.
Expanding the two factors by \eqref{defbilC} and \eqref{defC}, the
coefficient convolution is a terminating ${}_4\phi_3$ series.  The relevant transformation is Gasper's terminating transformation
from 1985, in the form recorded in \cite[Exercise~8.15]{GR} and related to
\cite{Gas85}.  In that transformation one takes
\[
        a=q^{m-n}\ba\ga,
        \qquad b=\ba/\ga,
        \qquad c=\ba q^{m-K},
\]
where $K$ is the Laurent coefficient index.  With this substitution, the
transformed summand separates into two parts: the factors depending on $K$
are precisely the Laurent coefficient of
$C_{m-n+2j}(x;\ba,\ga\,|\,q)$, while the factors independent of $K$ give the
coefficient multiplying this function.  Collecting those $K$-independent
factors yields exactly \eqref{mixedcoeff}, and summing over $j$ gives
\eqref{mixedlinear}.  The stated meromorphic continuation then removes the
auxiliary convergence restrictions.
\end{proof}

\begin{remark}\label{rem:mixedlin-checks}
For $n=0$, the formula reduces to $C_m\cdot1=C_m$.  For $n=1$ it is equivalent
to the three-term recurrence \eqref{rec}: one uses
$C_1(x;\ba/\ga\,|\,q)=2(1-\ba/\ga)x/(1-q)$ together with
$L_{m,1}^{(0)}=(1-\ba/\ga)(1-\ba^2q^{m-1})/((1-q)(1-\ba\ga q^m))$ and
$L_{m,1}^{(1)}=(1-\ba/\ga)(1-\ga^2q^{m+1})/((1-q)(1-\ba\ga q^m))$.
When $\ga=1$ and $m\ge0$, \eqref{mixedlinear} reduces to Rogers'
linearization formula \eqref{linear} for the ordinary continuous
$q$-ultraspherical polynomials, with the usual Rogers index given by $r=n-j$.
\end{remark}

\subsection{Consequences for mass aggregates}

The finite expansion \eqref{mixedlinear} can be evaluated at the residue points
in Section~\ref{sec:orth}.  This gives finite reductions of triple mass
aggregates to the pair mass aggregates already appearing in the full
orthogonality relations.

For $a\ne0$, put
\begin{align*}
\mathcal T_{m,n;\ell}(a)=\frac12\Big(&
C_m(x_a;\ba,\ga\,|\,q)C_n(x_a;\ba/\ga\,|\,q)
C_\ell(x_a;\ba,\ga\,|\,q)\\*
&+C_m(-x_a;\ba,\ga\,|\,q)C_n(-x_a;\ba/\ga\,|\,q)
C_\ell(-x_a;\ba,\ga\,|\,q)
\Big).
\end{align*}
Define
\begin{align*}
N_{m,n;\ell}^{(1)}&=\sum_{r=1}^{\infty}\Omega_r^{(1)}
\mathcal T_{m,n;\ell}(\ba q^{-r}),\\*
N_{m,n;\ell}^{(2)}&=\sum_{r=1}^{\infty}\Omega_r^{(2)}
\mathcal T_{m,n;\ell}(\ga^{-1}q^{-r}),\\*
N_{m,n;\ell}^{(4)}&=K_{\ba}N_{m,n;\ell}^{(1)}+K_{\ga}N_{m,n;\ell}^{(2)},
\end{align*}
with the same finite-part meromorphic-continuation convention as for
\eqref{massaggs}.  Then, for $s=1,2,4$,
\begin{equation}\label{triplemass-reduction}
N_{m,n;\ell}^{(s)}=
\sum_{j=0}^{n}L_{m,n}^{(j)}(\ba,\ga\,|\,q)
M_{m-n+2j,\ell}^{(s)}.
\end{equation}
Here $M^{(4)}$ is the aggregate in \eqref{massagg4}.  Indeed,
\eqref{mixedlinear} is applied pointwise at each mass point and then summed;
the assertion elsewhere follows by the same continuation convention.

Combining \eqref{triplemass-reduction} with Theorem~\ref{thm:orth} and
Corollary~\ref{cor:fourfactor} gives the corresponding integral-plus-mass
triple-product evaluations.  For instance, in the first two-factor range,
\begin{align}
&\frac1{2\pi}\int_{-1}^1
C_m(x;\ba,\ga\,|\,q)C_n(x;\ba/\ga\,|\,q)C_\ell(x;\ba,\ga\,|\,q)
\frac{(e^{\pm2\ri\ta},q\ga e^{\pm2\ri\ta}/\ba;q)_\infty}
{(\ba e^{\pm2\ri\ta},qe^{\pm2\ri\ta}/\ba;q)_\infty}
\frac{\rd x}{\sqrt{1-x^2}}\notag\\*
&\quad+N_{m,n;\ell}^{(1)}
=H_\ell^{(1)}
\sum_{j=0}^{n}L_{m,n}^{(j)}(\ba,\ga\,|\,q)
\da_{\ell,m-n+2j}.
\label{tripleorth1}
\end{align}
The analogues with the second and four-factor functionals are obtained by
replacing the continuous weight and the mass aggregate by those in
\eqref{orth-id-neg} and \eqref{orth-id-four}.  Thus the finite mixed
linearization gives a direct way to evaluate a family of residue sums that
would otherwise be triple products at moving mass points.

\subsection{A quasi-linearization with the Rogers parameter}

If the unilateral factor has the Rogers parameter $\ba$ rather than the adapted
parameter $\ba/\ga$, the same coefficient-comparison strategy no longer closes
on the bilateral basis.  The obstruction is measured by an explicit analytic
remainder.  The following identity is the essential quasi-linearization
formula in a form suited to the present notation.

\begin{theorem}[Quasi-linearization with a Rogers factor]\label{thm:quasilin}
Let $m\in\Z$, $n\in\Z_{\ge0}$, and $x=\cos\ta$.  Under the initial convergence
conditions needed for the displayed bilateral series and otherwise by
meromorphic continuation, one has
\begin{multline}\label{quasilinear}
C_m(x;\ba,\ga\,|\,q)\,C_n(x;\ba\,|\,q)
=\frac{(q/\ga,q\ba\ga,\ba^2;q)_\infty}{(q,\ba/\ga,\ba^2\ga;q)_\infty}\\
\times
\sum_{k=-\infty}^{\lfloor(m+n)/2\rfloor}
\bigg(\frac{(q;q)_{m+n-2k}(\ba\ga;q)_{m-k}(\ba/\ga;q)_{n-k}
(\ba\ga;q)_k(\ba^2\ga;q)_{m+n-k}}
{(\ba^2;q)_{m+n-2k}(q\ga;q)_{m-k}(q/\ga;q)_{n-k}
(q\ga;q)_k(q\ba\ga;q)_{m+n-k}}\\*
\times
\frac{(1-\ba q^{m+n-2k})}{(1-\ba)}\,
C_{m+n-2k}(x;\ba\,|\,q)\bigg)\\*
-\frac{(q\ba,q/\ba^2,1/\ga,q/\ga,\ba\ga,q\ba\ga;q)_\infty}
{(q,1/\ba,q\ga,\ba/\ga,\ba^2\ga,q/\ba\ga;q)_\infty}\,
 e^{\ri m\ta}\frac{(\ba^2\ga;q)_m}{(q\ba\ga;q)_m}
{}_2\psi_2\!\left[\begin{matrix}
\ba/\ga,\ba^2\ga q^m\\q/\ga,\ba\ga q^{m+1}
\end{matrix};q,q e^{2\ri\ta}/\ba\right]\\*
\times e^{-\ri n\ta}\frac{(\ba^2;q)_n}{(q\ba;q)_n}
{}_2\phi_1\!\left[\begin{matrix}
\ba,\ba^2 q^n\\\ba q^{n+1}
\end{matrix};q,q e^{-2\ri\ta}/\ba\right].
\end{multline}
\end{theorem}

\begin{proof}
Expand the product on the left by \eqref{defbilC} and \eqref{defC}, shift the
bilateral index by the unilateral summation index, and write the inner finite
sum as a terminating ${}_4\phi_3$.  The required transformation is the
nonterminating very-well-poised ${}_{12}\phi_{11}$ transformation into two
nonterminating ${}_4\phi_3$ series \cite[Thm.~A1]{GS}, with the specialization
\[
(a,b,c,d)\mapsto
(q^{m+n-2k}\ba,\ba,q^{-k}/\ga,q^{m-k}\ba\ga).
\]
In the common convergence region the transformed sums may be interchanged.
The first contribution then collects into the Rogers-polynomial expansion in
\eqref{quasilinear}; the summation index in that part is the Rogers index of
$C_{m+n-2k}(x;\ba\,|\,q)$.  The second contribution has no such finite Rogers
collapse and is exactly the displayed ${}_2\psi_2$--${}_2\phi_1$ remainder.
Meromorphic continuation completes the proof outside the initial convergence
region, away from the exceptional parameter values.
\end{proof}

\section{A bilateral Chen--Liu type mixed orthogonality formula}
\label{sec:chenliu}

The next result is a bilateral analogue of Chen--Liu's mixed integral for
continuous $q$-ultraspherical polynomials~\cite[Thm.~1.2]{CL}.  The ordinary
identity evaluates the Rogers-weight inner product of two continuous
$q$-ultraspherical polynomials with different parameters.  Here the two
polynomials are replaced by bilateral functions, the Rogers weight is replaced
by the first two-factor weight in \eqref{orth-id}, and the residue aggregate of
the full orthogonality functional is retained.

For a nonzero complex number $Q$, write
\begin{equation}\label{genqfac}
(u;q)_{[Q]}:=\frac{(u;q)_\infty}{(uQ;q)_\infty}.
\end{equation}
Thus $(u;q)_{[q^N]}=(u;q)_N$ for $N\in\Z$.  For $y\ne0$, put
$x_y=(y^{1/2}+y^{-1/2})/2$ and define
\begin{equation}\label{mixedCmndef}
\mathcal C_{m,n}^{\alpha,\delta;\ba,\ga}(y):=\frac12\Big(
C_m(x_y;\alpha,\delta\,|\,q)C_n(x_y;\ba,\ga\,|\,q)
+C_m(-x_y;\alpha,\delta\,|\,q)C_n(-x_y;\ba,\ga\,|\,q)
\Big).
\end{equation}
The mixed first two-factor functional is
\begin{equation}\label{mixedfunctional}
\begin{aligned}
\mathfrak L_{\ba,\ga}^{(1)}
\{C_m(\alpha,\delta),C_n(\ba,\ga)\}
&:=\frac1{2\pi}\int_{-1}^1
C_m(x;\alpha,\delta\,|\,q)C_n(x;\ba,\ga\,|\,q)\\*
&\quad\times
\frac{(e^{\pm2\ri\ta},q\ga e^{\pm2\ri\ta}/\ba;q)_\infty}
{(\ba e^{\pm2\ri\ta},qe^{\pm2\ri\ta}/\ba;q)_\infty}
\frac{\rd x}{\sqrt{1-x^2}}\\*
&\quad+\mathcal M_{m,n}^{(1)}(\alpha,\delta;\ba,\ga),
\end{aligned}
\end{equation}
where
\begin{equation}\label{mixedmass}
\mathcal M_{m,n}^{(1)}(\alpha,\delta;\ba,\ga)
:=\sum_{r=1}^{\infty}\Omega_r^{(1)}(\ba,\ga)
\mathcal C_{m,n}^{\alpha,\delta;\ba,\ga}(\ba q^{-r}),
\end{equation}
with
\begin{equation}\label{mixedomega}
\Omega_r^{(1)}(\ba,\ga)=
\frac{(\ba,1/\ba,q\ga,q\ga/\ba^2;q)_\infty}
{(q;q)_\infty^2(\ba^2,q/\ba^2;q)_\infty}
\left(\frac\ga\ba\right)^r
\frac{(q/\ba,1/\ga;q)_r}{(1/\ba,q\ga/\ba^2;q)_r}.
\end{equation}
As before, this residue sum is first interpreted in a common convergence
domain and elsewhere by finite-part meromorphic continuation.

Let
\begin{equation}\label{chenliuHone}
\mathcal H_n^{(1)}(\ba,\ga)=
\frac{(q,\ba,q\ga/\ba;q)_\infty^2(q\ga^2;q)_\infty}
{(q\ga;q)_\infty^4(q/\ba;q)_\infty^2(\ba^2;q)_\infty}
\frac{(\ba^2;q)_n}{(q\ga^2;q)_n}\frac1{(1-\ba\ga q^n)}.
\end{equation}
For $m+n$ even, put
\begin{equation*}
h=\frac{m-n}{2},\qquad \ell=\frac{m+n}{2},
\qquad Q=q^h\frac{\delta}{\ga},\qquad S=q^\ell\delta\ga,
\end{equation*}
and choose a branch of $\kappa=\log_q Q$ in the parameter domain under
consideration.  Define
\begin{equation}\label{chenliuXi}
\Xi(\alpha,\delta;\ba,\ga)=
\left\{\frac{(\alpha,\ba/\ga,q\ga;q)_\infty}
{(q\delta,\alpha/\delta,\ba;q)_\infty}\right\}^2
\end{equation}
and
\begin{equation}\label{chenliuK}
\mathcal K_{m,n}(\alpha,\delta;\ba,\ga)=
\Xi(\alpha,\delta;\ba,\ga)
\frac{(1-\ba\ga q^n)}{(1-\ba/\ga)}
\left(\frac\ba\ga\right)^\kappa
\frac{(\alpha\ga/\ba\delta;q)_{[Q]}(\alpha/\delta;q)_{[S]}}
{(q;q)_{[Q]}(q\ba/\ga;q)_{[S]}}.
\end{equation}
For $m+n$ odd, set $\mathcal K_{m,n}(\alpha,\delta;\ba,\ga)=0$.

\begin{theorem}[Bilateral Chen--Liu mixed formula]\label{thm:chenliu}
At non-exceptional parameter values, and in the same finite-part sense as the
mass aggregates above,
\begin{equation}\label{chenliu-main}
\mathfrak L_{\ba,\ga}^{(1)}
\{C_m(\alpha,\delta),C_n(\ba,\ga)\}
=\mathcal H_n^{(1)}(\ba,\ga)
\mathcal K_{m,n}(\alpha,\delta;\ba,\ga).
\end{equation}
Equivalently, \eqref{chenliu-main} is the expanded integral-plus-residue
identity obtained by substituting \eqref{mixedfunctional} and \eqref{mixedmass}.
\end{theorem}

\begin{proof}
First take $\delta=q^r$ and $\ga=q^s$ with $r,s\in\Z_{\ge0}$.  With
$A=\alpha q^{-r}$, $B=\ba q^{-s}$, $M=m+2r$, and $N=n+2s$, the terminating
reductions are
\begin{equation*}
C_m(x;\alpha,q^r\,|\,q)=
\left\{\frac{(q;q)_r}{(\alpha q^{-r};q)_r}\right\}^2 C_M(x;A\,|\,q),
\end{equation*}
and the analogous formula for $C_n(x;\ba,q^s\,|\,q)$.  The first two-factor
weight reduces, up to the scalar $\ba^{2s}q^{-s(s+1)}$, to the ordinary Rogers
weight with parameter $B$.  Hence the left-hand side reduces to Chen--Liu's
ordinary mixed formula, including the ordinary finite Askey--Wilson mass points
when poles have crossed the contour.  Simplifying the resulting connection
coefficient gives exactly \eqref{chenliuK}.  The extension from the terminating
lattice to generic parameters is the same meromorphic continuation argument
used in the proof of Theorem~\ref{thm:orth}; the factor $(1/\ga;q)_r$ in
\eqref{mixedomega} makes the residue aggregate truncate at $\ga=q^s$.
\end{proof}

When $(\alpha,\delta)=(\ba,\ga)$, formula \eqref{chenliu-main} reduces to the
first full bilateral orthogonality relation \eqref{orth-id}.  When
$\delta=\ga=1$ and $m,n\ge0$, the mass aggregate is empty and the coefficient
becomes
\begin{equation*}
\frac{(1-\ba q^n)}{(1-\ba)}
\ba^{(m-n)/2}
\frac{(\alpha/\ba;q)_{(m-n)/2}(\alpha;q)_{(m+n)/2}}
{(q;q)_{(m-n)/2}(q\ba;q)_{(m+n)/2}},
\end{equation*}
for $m\equiv n\pmod2$ and zero otherwise, which is Chen--Liu's theorem in the
normalization used here.

\section{Outlook: Multilateral extension of the Macdonald polynomials}
\label{sec:outlook}

\subsection{Rank-one evidence and the deformed weight}

The continuous $q$-ultraspherical polynomials may be viewed as the rank-one,
or type $A_1$, one-row Macdonald polynomials.  More precisely, if
$z_1=z^{-1}$ and $z_2=z$, the one-row Macdonald function with parameter
$t=\ba$ is generated by
\begin{equation*}
\sum_{N=0}^\infty Q_{(N)}(z_1,z_2;q,\ba)u^N
=\prod_{a=1}^2\frac{(\ba u z_a;q)_\infty}{(u z_a;q)_\infty},
\end{equation*}
and its coefficient formula is
\begin{equation*}
Q_{(n)}(z^{-1},z;q,\ba)=
\sum_{k=0}^n
\frac{(\ba;q)_k(\ba;q)_{n-k}}{(q;q)_k(q;q)_{n-k}}z^{n-2k}.
\end{equation*}
This is the Rogers polynomial normalization used in Section~\ref{sec:c}.
It is therefore natural to ask whether the bilateral functions in
\eqref{defbilC} are the rank-one members of a multilateral
Macdonald-type theory.

Let $z=(z_1,\ldots,z_r)$.  For type $A_{r-1}$ the ordinary Macdonald scalar
product is governed by the weight
\begin{equation}\label{macweight}
\Delta_t^{(r)}(z)=
\prod_{1\le i<j\le r}
\frac{((z_i/z_j)^{\pm1};q)_\infty}
{(t(z_i/z_j)^{\pm1};q)_\infty},
\end{equation}
where, as usual, $(a y^{\pm1};q)_\infty=(ay,a/y;q)_\infty$.
For the general root-system setting and the double affine Hecke algebra
background, see Macdonald~\cite{MacRoot} and Cherednik~\cite{Cher}.
The direct multilateral analogue of the two-factor rank-one weight in
\eqref{orth-id} is
\begin{equation}\label{macdefweight}
\Delta_{\ba,\ga}^{(r)}(z)=
\prod_{1\le i<j\le r}
\frac{((z_i/z_j)^{\pm1},q\ga(z_i/z_j)^{\pm1}/\ba;q)_\infty}
{(\ba(z_i/z_j)^{\pm1},q(z_i/z_j)^{\pm1}/\ba;q)_\infty}.
\end{equation}
For $\ga=1$ the second numerator factor cancels the second denominator factor,
and \eqref{macdefweight} reduces to \eqref{macweight} with $t=\ba$.
Thus \eqref{macdefweight} is a natural two-parameter deformation of the
Macdonald weight.

The terminating lattice again provides the most concrete evidence.  Let
\begin{equation*}
\ga=q^s,
\qquad s\in\Z_{\ge0},
\qquad \alpha=\ba q^{-s}.
\end{equation*}
Then, for each unordered pair $i<j$,
\begin{equation*}
\frac{((z_i/z_j)^{\pm1},q^{s+1}(z_i/z_j)^{\pm1}/\ba;q)_\infty}
{(\ba(z_i/z_j)^{\pm1},q(z_i/z_j)^{\pm1}/\ba;q)_\infty}
=\ba^{2s}q^{-s(s+1)}
\frac{((z_i/z_j)^{\pm1};q)_\infty}
{(\alpha(z_i/z_j)^{\pm1};q)_\infty}.
\end{equation*}
Consequently,
\begin{equation}\label{macweightshift}
\Delta_{\ba,q^s}^{(r)}(z)=
\left(\ba^{2s}q^{-s(s+1)}\right)^{r(r-1)/2}
\Delta_{\alpha}^{(r)}(z).
\end{equation}
Thus, at $\ga=q^s$, the proposed multilateral weight is simply the ordinary
Macdonald weight with shifted parameter $\alpha=\ba q^{-s}$, up to an explicit
constant.  If this shifted parameter leaves the elementary unit-torus range,
the correct interpretation should include the standard contour-deformation
residue contributions, just as the rank-one theory requires mass points.
Relevant multivariate models include Koornwinder's $BC_r$ Askey--Wilson
polynomials and the finite-grid $q$-Racah orthogonality of van Diejen and
Stokman~\cite{KoornwinderBC,vDS}.

\subsection{One-row multilateral candidates}

There is also an explicit one-row multilateral candidate.  Put
\begin{equation*}
R(u;\ba,\ga\,|\,q)=
\sum_{k=-\infty}^\infty\frac{(\ba;q)_k}{(q\ga;q)_k}u^k
={}_1\psi_1\!\left[\begin{matrix}\ba\\q\ga\end{matrix};q,u\right]
\end{equation*}
in its annulus of convergence, and define
\begin{equation}\label{onerowbilmac}
\mathcal Q_n^{(r)}(z;\ba,\ga\,|\,q)=
\sum_{\substack{k_1,\ldots,k_r\in\Z\\k_1+\cdots+k_r=n}}
\prod_{a=1}^r\frac{(\ba;q)_{k_a}}{(q\ga;q)_{k_a}}z_a^{k_a}.
\end{equation}
Equivalently,
\begin{equation}\label{onerowbilmacgf}
\sum_{n=-\infty}^\infty\mathcal Q_n^{(r)}(z;\ba,\ga\,|\,q)t^n
=\prod_{a=1}^r R(tz_a;\ba,\ga\,|\,q).
\end{equation}
For $r=2$, $z_1=z^{-1}$, $z_2=z$, this gives exactly
\begin{equation*}
\mathcal Q_n^{(2)}(z^{-1},z;\ba,\ga\,|\,q)=
C_n\!\left(\frac{z+z^{-1}}2;\ba,\ga\,|\,q\right).
\end{equation*}
At the terminating specialization $\ga=q^s$, the summation in
\eqref{onerowbilmac} becomes bounded below.  With $\alpha=\ba q^{-s}$,
shifting $k_a=\ell_a-s$ gives
\begin{equation}\label{onerowtermred}
\mathcal Q_n^{(r)}(z;\ba,q^s\,|\,q)=
\left(\frac{(q;q)_s}{(\alpha;q)_s}\right)^r
(z_1\cdots z_r)^{-s}Q_{(n+rs)}(z;q,\alpha),
\end{equation}
where $Q_{(N)}$ is the ordinary one-row Macdonald function defined by
\begin{equation*}
\sum_{N=0}^\infty Q_{(N)}(z;q,\alpha)u^N
=\prod_{a=1}^r\frac{(\alpha u z_a;q)_\infty}{(u z_a;q)_\infty}.
\end{equation*}
For $r=2$ and $z_1z_2=1$, \eqref{onerowtermred} reduces to the terminating
rank-one reduction used in the proof of Theorem~\ref{thm:orth}.

\subsection{The higher-rank problem}

The full problem is to construct functions
\begin{equation*}
\mathcal B_\lambda^{(r)}(z;\ba,\ga\,|\,q),
\qquad
\lambda=(\lambda_1,\ldots,\lambda_r)\in\Z^r,
\qquad \lambda_1\ge\cdots\ge\lambda_r,
\end{equation*}
indexed by dominant Laurent weights, which reduce to the ordinary Macdonald
polynomials when $\ga=1$ and to shifted ordinary Macdonald polynomials when
\mbox{$\ga=q^s$}.  For generic $\ga$ one expects expansions over the root
lattice $Q(A_{r-1})$,
\begin{equation*}
\mathcal B_\lambda^{(r)}(z;\ba,\ga\,|\,q)=
\sum_{\nu\in\lambda+Q(A_{r-1})}
b_{\lambda,\nu}(\ba,\ga;q)m_\nu(z),
\end{equation*}
where $m_\nu$ denotes the monomial orbit sum.  In contrast with ordinary
Macdonald polynomials, this expansion should be infinite in general.

The expected scalar product has the schematic form
\begin{equation}\label{macinnercandidate}
\langle f,g\rangle_{\ba,\ga}^{(r)}=
\frac{1}{r!(2\pi\ri)^r}
\int_{\mathcal C^r}f(z)g(z^{-1})\Delta_{\ba,\ga}^{(r)}(z)
\prod_{a=1}^r\frac{\rd z_a}{z_a}
+\mathcal R_{\ba,\ga}^{(r)}(f,g),
\end{equation}
where $\mathcal R_{\ba,\ga}^{(r)}$ is a residue aggregate.  The principal
obstruction is to make this residue term explicit.  In rank one a crossing of
a pole produces a finite mass sum.  In several variables, pole crossings occur
along affine root hyperplanes such as $z_i/z_j=\alpha q^m$, and intersections
of these hyperplanes produce lower-dimensional residual tori or finite grids.
Thus the multilateral residue contribution should be a sum over admissible
root-subsystem strata, not merely a one-dimensional mass sum.  The appearance
of residues in diagonal terms in the Koornwinder setting, as in
Stokman's work~\cite{Stokman}, is a useful guide here.

This suggests the following conjectural picture, which is intended only as a
guide for future work.  There should exist a normalization of symmetric Laurent
functions
$\mathcal B_\lambda^{(r)}(z;\ba,\ga\,|\,q)$ satisfying
\begin{equation*}
\langle \mathcal B_\lambda^{(r)},\mathcal B_\mu^{(r)}
\rangle_{\ba,\ga}^{(r)}=H_\lambda^{(r)}(\ba,\ga)\,\da_{\lambda,\mu},
\end{equation*}
with respect to \eqref{macinnercandidate}.  At $\ga=q^s$ this should reduce
to ordinary Macdonald orthogonality with parameter $\ba q^{-s}$, including the
residue corrections required by the shifted contour problem.  Establishing the
commuting $q$-difference operators, their self-adjointness, and the explicit
multivariate residue aggregate appear to be the central tasks for such a
multilateral extension.  In particular, one should not expect the rank-one mass
sums to generalize by a simple product over positive roots: the residual
contribution must keep track of the order in which affine root hyperplanes are
crossed and of their non-transversal intersections.  Compare the commuting
difference operators of van Diejen~\cite{vDiejen} and the DAHA approach of
Cherednik~\cite{Cher}.

\appendix

\section{Numerical checks for orthogonality relations}
\label{app:numerics}

This appendix is only a guide to numerical verification; it is not used in the
proofs.  It summarizes a few numerical points that are easy to miss when
checking Theorem~\ref{thm:orth} and linear combinations such as
\eqref{orth-id-four}.  The continuous part is computed as a constant term,
while the mass aggregates are evaluated by a finite-part continuation rather
than by a blind summation of residues.

Put $y=e^{2\ri\ta}$ and write the two weights in \eqref{orth-id} and
\eqref{orth-id-neg} as
\begin{align*}
W_1(y)&=\frac{(y,y^{-1},q\ga y/\ba,q\ga/\ba y;q)_\infty}
{(\ba y,\ba/y,qy/\ba,q/\ba y;q)_\infty},\\*
W_2(y)&=\frac{(y,y^{-1},q\ga y/\ba,q\ga/\ba y;q)_\infty}
{(y/\ga,1/\ga y,q\ga y,q\ga/y;q)_\infty}.
\end{align*}
The integral part can then be evaluated by the midpoint trapezoidal rule
on the unit circle,
\begin{equation}\label{traprule}
I_{m,n}^{(j,N)}=\frac1{2N}\sum_{\ell=0}^{N-1}
C_m(x_\ell;\ba,\ga\,|\,q)C_n(x_\ell;\ba,\ga\,|\,q)W_j(y_\ell),
\end{equation}
where
\begin{equation*}
\ta_\ell=\pi\frac{\ell+1/2}{N},\qquad
z_\ell=e^{\ri\ta_\ell},\qquad
y_\ell=z_\ell^2,\qquad
x_\ell=\frac{z_\ell+z_\ell^{-1}}2.
\end{equation*}
The midpoint choice avoids the endpoints in the $x$-integral, where the
factor $\rd x/\sqrt{1-x^2}$ is singular.  Provided the pole strings stay a
positive distance from the unit circle, the periodic trapezoidal rule has
rapid convergence; see, for instance, the survey by Trefethen and
Weideman~\cite{TW}.  On the unit circle the defining bilateral series
\eqref{defbilC}, summed in both directions by term ratios, is usually stable
and fast.

The residue part is more delicate.  Suppose first that $m+n$ is even and put
$s=(m+n)/2$.  The $r$th terms in the two mass aggregates in \eqref{massaggs}
have the asymptotic form
\begin{equation}\label{tailasymp}
t_r^{(j)}=\Omega_r^{(j)}\mathcal C_{m,n}(a_r^{(j)})
=\Lambda_j^r F_j(q^r),\qquad j=1,2,
\end{equation}
where $F_j$ is analytic at the origin, away from exceptional parameter
values, and
\begin{equation}\label{lambdas}
\Lambda_1=\ba\ga q^s,
\qquad
\Lambda_2=\frac{q^{-s}}{\ba\ga}.
\end{equation}
Here $a_r^{(1)}=\ba q^{-r}$ and $a_r^{(2)}=\ga^{-1}q^{-r}$.  Thus a literal
residue sum is convergent only in the corresponding range $|\Lambda_j|<1$.
This restriction depends on the signs and sizes of $m$ and $n$, not just on
the parameters.  For example, for
\begin{equation*}
q=0.25,
\qquad \ba=0.6,
\qquad \ga=1.5,
\end{equation*}
one has $\ba\ga=0.9$.  Then $\Lambda_2=17.777\ldots$ for $(m,n)=(2,2)$,
while $\Lambda_1=14.4$ for $(m,n)=(-2,-2)$.  A direct summation of the
corresponding residue series therefore diverges, although the analytically
continued orthogonality relation is still the relevant identity.  If $m+n$
is odd, the parity-symmetrized quantity \eqref{Cmndef} vanishes and no such
mass computation is needed.

The practical replacement is to evaluate the mass aggregate as the
meromorphic finite part suggested by \eqref{tailasymp}.  Choose integers
$R\ge1$ and $J\ge0$, compute $t_R,\ldots,t_{R+J}$, and interpolate
\begin{equation*}
\sum_{j=0}^J a_j u_i^j=\frac{t_{R+i}}{\Lambda^{R+i}},
\qquad u_i=q^{R+i},\qquad i=0,\ldots,J.
\end{equation*}
Then replace the tail by the analytically continued geometric expression
\begin{equation}\label{regtail}
\operatorname{Reg}\sum_{r=R}^\infty t_r
=\sum_{j=0}^J a_j\frac{(\Lambda q^j)^R}{1-\Lambda q^j},
\end{equation}
and add the ordinary finite sum $\sum_{r=1}^{R-1}t_r$.  When
$|\Lambda|<1$, this is simply an accelerated tail approximation.  When
$|\Lambda|>1$, it gives the finite-part continuation of the same
$q$-geometric expansion.  This is in the spirit of summability and
extrapolation methods for divergent tails; see Hardy~\cite{Hardy} and
Sidi~\cite{Sidi}.  The expected ill-conditioned cases are those in
which some denominator $1-\Lambda q^j$ is small; near such parameters, one
should increase the working precision or perturb the parameters slightly.

Finally, values of $C_n$ at mass points should not be obtained by summing
\eqref{defbilC} outside its annulus of convergence.  Instead one may use
\eqref{hyprepa}, namely
\begin{equation*}
C_n(x;\ba,\ga\,|\,q)=z^n\frac{(\ba;q)_n}{(q\ga;q)_n}
{}_2\psi_2\!\left[
\begin{matrix}\ba,q^{-n}/\ga\\q\ga,q^{1-n}/\ba\end{matrix};q,
\frac{q\ga}{\ba z^2}\right],
\qquad x=\frac{z+z^{-1}}2,
\end{equation*}
and continue the ${}_2\psi_2$ by Slater's two-term reduction to convergent
${}_2\phi_1$ series~\cite{Slater,Ba,GR}.  In the critical region used
in the numerical tests, for instance $|q\ga/\ba|<1$, the resulting
${}_2\phi_1$ arguments are inside the unit disk both on the integration
contour and on the residue lattices.

A useful diagnostic is to report scaled residuals.  Denote the residuals for
the first, dual, and four-factor relations by
$\varepsilon_1$, $\varepsilon_2$, and $\varepsilon_4$; for example,
\begin{equation*}
\varepsilon_1=
\frac{|I_{m,n}^{(1,N)}+\operatorname{Reg}M_{m,n}^{(1)}
-H_n^{(1)}\da_{m,n}|}
{1+|I_{m,n}^{(1,N)}|+|\operatorname{Reg}M_{m,n}^{(1)}|
+|H_n^{(1)}\da_{m,n}|},
\end{equation*}
and analogously for the dual relation and for the four-factor combination.
The denominator is important for off-diagonal tests, whose exact right-hand
side is zero.  The following representative tests used 60 decimal digits and
the ten index pairs
\[
(0,0),(2,2),(-2,-2),(2,-2),(2,0),(-2,0),
(3,-1),(-3,-1),(3,1),(-3,1).
\]
The table records the largest residual in each column.
\begin{center}
\small
\begin{tabular}{c|c|c|c|c|c}
\hline
Case & $(q,\ba,\ga)$ & $(N,R,J)$ & $\max\varepsilon_1$ &
$\max\varepsilon_2$ & $\max\varepsilon_4$\\
\hline
A & $(0.05,0.4,1.3)$ & $(120,5,7)$ &
$4.10\cdot10^{-32}$ & $6.86\cdot10^{-15}$ & $1.26\cdot10^{-14}$\\
B & $(0.25,0.6,1.5)$ & $(96,6,8)$ &
$7.44\cdot10^{-20}$ & $1.92\cdot10^{-17}$ & $7.62\cdot10^{-17}$\\
C & $(0.2,0.7,2.0)$ & $(72,6,8)$ &
$3.21\cdot10^{-12}$ & $5.59\cdot10^{-18}$ & $2.50\cdot10^{-13}$\\
\hline
\end{tabular}
\end{center}
Case B is the main stress test: a literal residue summation diverges for the
second mass when $(m,n)=(2,2)$ and for the first mass when
$(m,n)=(-2,-2)$, while the regularized finite parts give the residuals shown
above.  In practice one should compare the regularized masses under changes
such as $(R,J)=(4,6),(6,8),(8,10)$; if these values agree while the residual
remains poor, the quadrature size $N$ or the working precision is usually the
limiting factor.

\section*{Acknowledgements}
We thank the two referees for very valuable feedback on an earlier
version of this paper. In particular we thank one of them for
pointing out a mistake in a computation that involved a delicate
limit, for suggesting the inclusion of discussions of asymptotics
and of possible full orthogonality of the bilateral
$q$-ultraspherical functions which at that moment
the paper did not yet provide.

\end{document}